\documentclass[smallcondensed]{svjour3}
\usepackage{graphicx}
\usepackage{caption,color}
\usepackage{amsfonts}
\usepackage{subfigure}
\usepackage{mathrsfs}
\usepackage{amssymb,amsmath}

\numberwithin{equation}{section}

\newcommand{\p}{\partial}
\newcommand{\og}{\omega}
\newcommand{\sech}{{\mathrm{sech}} }
\newcommand{\Og}{\Omega}
\newcommand{\fl}[2]{\frac{#1}{#2}}
\newcommand{\be}{\begin{equation}}
\newcommand{\ee}{\end{equation}}
\newcommand{\nn}{\nonumber}
\newcommand{\wh}{\widehat}
\newcommand{\wt}{\widetilde}

\newcommand{\tht}{\theta}

\renewcommand{\thefigure}{\thesection.\arabic{figure}}
\renewcommand{\thetable}{\thesection.\arabic{table}}

\begin{document}

\title{A Deuflhard-type Exponential Integrator Fourier Pseudospectral Method for the ``Good" Boussinesq Equation\thanks{This first author was supported by the Alexander von Humboldt Foundation and the second author was supported by NSFC (No. 11801183).}
}

\titlerunning{An exponential integrator for the GB equation}

\author{Chunmei Su \and
Wenqi Yao}

\institute{ Chunmei Su\at
Zentrum Mathematik, Technische Universit\"{a}t M\"unchen, 85748 Garching bei M\"unchen, Germany\\
              \email{sucm13@163.com}
   \and
Wenqi Yao\at
School of Mathematics, South China University of Technology, 510641 Guangzhou, China\\
\email{yaowq@scut.edu.cn}
}
\date{}
\maketitle

\begin{abstract}
We propose an exponential integrator Fourier pseudospectral method DEI-FP for solving the ``Good" Boussinesq (GB)  equation. The numerical scheme is based on a Deuflhard-type exponential integrator and a Fourier pseudospectral method for temporal and spatial discretizations, respectively. The scheme is fully explicit and efficient due to the fast Fourier transform. Rigorous error estimates are established for the  method without any CFL-type condition constraint. In more details, the method converges quadratically and spectrally in time and space, respectively. Extensive numerical experiments are reported to confirm the theoretical analysis and to demonstrate rich dynamics of the GB equation.

\keywords{Error estimate \and Exponential integrator \and ``Good" Boussinesq equation  \and Pseudospectral method}

\subclass{35Q53 \and 65M15 \and 65M70}
\end{abstract}

\section{Introduction}
Similar to the Korteweg--de Vries (KdV) equation or cubic Schr\"odinger equation, the ``Good" Boussinesq (GB) equation \cite{boussinesq1872} also gives rise to solitons:
\be\label{GB}
z_{tt}-z_{xx}+z_{xxxx}-(z^2)_{xx}=0.
\ee
However, the GB equation possesses some remarkable features which make it different from the KdV equation, e.g., two solitary waves can merge into a single wave or develop into the so-called antisolitons \cite{manoranjan1988}. The GB equation was originally introduced to model one-dimensional weakly nonlinear dispersive waves in shallow water \cite{boussinesq1872}. Furthermore, it was extended by replacing the quadratic nonlinearity with a general function of $z$ to model small oscillations of nonlinear beams \cite{Varl} or the two-way propagation of water waves in a channel \cite{Bona}.

%As is well-known, th GB equation conserves the energy
%\[E:=\int (w^2+z^2+z_x^2+\fl{2}{3}z^3)dx,\quad w_x=z_t,\]
%for the whole-space problem or finite interval problem under the periodic or homogeneous Dirichlet boundary conditions \cite{Chen2017,yan2016}.

There have been extensive studies for the GB equation in recent decades. For the local well-posedness of the initial value problem, we refer to \cite{Bona1988,farah2009,Kish2010}. Specifically, the GB equation is locally well posed for initial data in $H^s(\mathbb{R})\times H^{s-2}(\mathbb{R})$ with $s>-1/2$.
In the periodic setting, we refer to \cite{fang1996,farah2010,kishimoto2013sharp,oh2013improved} for the local well-posedness. More precisely, it was shown in \cite{kishimoto2013sharp} that the initial value problem in $H^s(\mathbb{T})\times H^{s-2}(\mathbb{T})$ is locally well-posed for $s\ge -1/2$ and ill-posed for $s<-1/2$. For the numerical and analytical study of the GB equation, it can be traced back to \cite{manoranjan1984}, where the soliton interaction mechanism was investigated and some numerical experiments were reported while little analysis was given on the stability and convergence of the methods employed. Since then numerous numerical methods were developed for solving the GB equation.
 Finite difference methods (FDM) have been developed in \cite{bratsos2007,Zoh,ortega1990}. Particularly, in \cite{ortega1990}, the nonlinear stability and convergence was shown for a family of explicit finite difference schemes for solving the GB equation under a severe CFL-type condition: $\Delta t=O(\Delta x^2)$ where $\Delta t$ and $\Delta x$ represent the discretization parameters in time and space, respectively. More competitive Fourier spectral methods were also proposed and analyzed for the GB equation \cite{cheng2015,de1990,de1991,yan2016}. A second order temporal pseudospectral discretization was proposed in \cite{de1991} and the full order convergence was proved in a weak energy norm: the $L^2$ norm in $z$ combined with the $H^{-2}$ norm in $z_t$ under a similar time step constraint $\Delta t=O(\Delta x^2)$.  The energy norm was improved in \cite{cheng2015}, where a second order temporal scheme was proved to converge unconditionally in a stronger energy norm: the $H^2$ norm in $z$ combined with the $L^2$ norm in $z_t$. However, it requires the solution to be regular enough: $z\in H^4([0,T];L^2)\cap L^\infty([0,T];H^{m+4})\cap H^2([0,T];H^4)$ to get an error of $O(\Delta t^2+\Delta x^m)$. There are some other methods in the literatures for solving the GB equation, such as Petrov-Galerkin methods \cite{manoranjan1984}, meshless methods \cite{dehghan2012}, operator splitting methods \cite{zhang2018,zhang2017}, energy-preserving methods \cite{cai2013,Chen2017,jiang2016,yan2016} and Runge-Kutta exponential integrators \cite{Akbar}. In more details, an unconditional full order convergence was shown for a second order operator splitting numerical scheme in the energy norm \cite{zhang2017}: the $H^2 \times L^2$ norm in $z \times z_t$.

Nowadays, exponential time integrators have been widely applied for parabolic and hyperbolic equations \cite{exp10,Akbar,Alex19,zhao2016error}. Particularly, several efficient schemes were proposed for solving the GB equation \cite{Akbar} based on fourth-order exponential integrators of Runge-Kutta type. However, nothing concerning the stability and convergence results was involved. Very recently, two low regularity exponential-type integrators were proposed and analyzed for the GB equation \cite{Alex19} based on a technique of \emph{twisting variable}. The advantage of the low regularity exponential-type integrators is that it requires less regularities on the solution to obtain the same accuracy, compared to the classical exponential integrators. However, the design of low regularity integrators depends strongly on the particular form of the nonlinearity. Such an integrator can hardly be extended to more general equations, e.g., the GB equation with a general nonlinearity.

In the present work, we consider a second-order Deuflhard-type exponential integrator for solving the GB equation. Our work differs from the existing studies in the literatures in the following two aspects: (i) we consider the GB equation with a general nonlinearity; (ii) an unconditional convergence is proved in a general energy norm: the $H^m$ norm in $z$ combined with the $H^{m-2}$ norm in $z_t$ for $m>1/2$. The rest of this paper is organized as follows. In Section 2, we propose the Deuflhard-type exponential integrator Fourier pseudospectral (DEI-FP) method for the GB equation. The main error estimate is given and proved in Section 3. Numerical results are reported in Section 4 to illustrate the proved convergence results and to demonstrate the complicated dynamics of the GB equation. Finally, some concluding remarks are drawn in Section 5.

\section{A Deuflhard-type exponential integrator Fourier pseudospectral method}
In this section, we present the exponential integrator Fourier pseudospectral (DEI-FP)
method for the GB equation, based on a Deuflhard-type time integrator in combination with a Fourier pseudospectral discretization in space.

For implementation issues, we consider the GB equation with a general nonlinearity with periodic boundary conditions imposed:
\be \label{gB}
\left\{
\begin{aligned}
&z_{tt}-z_{xx}+z_{xxxx}-(f(z))_{xx}=0,  \quad x \in \Og,\quad t>0,\\
&z(x,0)=z_0(x),\quad \partial_t z(x,0)=z_1(x),\quad x\in \overline{\Omega},
 \end{aligned}
         \right .
\ee
where $f(\cdot)$ is a smooth functional. For integer $m>0$,
$\Omega=[a, b]$, we denote by $H^m(\Omega)$  the standard Sobolev space with norm
\be\label{sn}
\|f\|_m^2=\sum\limits_{l} (1+|\mu_l|^2)^m|\widehat{f}_l|^2,\quad \mathrm{for}\quad f(x)=\sum\limits_{l\in \mathbb{Z}} \widehat{f}_l e^{i\mu_l(x-a)},\quad \mu_l=\fl{2l\pi }{b-a}.
\ee
For $m=0$, the space is exactly $L^2(\Og)$ and the corresponding norm is denoted as $\|\cdot\|$.
Furthermore, we denote by $H_{\mathrm p}^m(\Omega)$ the subspace of $H^m(\Omega)$ which consists of functions with derivatives of order up to $m-1$ being $(b-a)$-periodic.
We see that the space $H^m(\Omega)$ with fractional $m$ is also well-defined which consists of functions such that $\|\cdot\|_m$ is finite.

For the full discretization of \eqref{GB}, we introduce some discrete spaces. Choose a mesh size $h:=(b-a)/M$ with $M$ a positive integer, and a time step $\tau>0$. Denote the grid points and time steps as
$$x_j:=a+jh,\quad j=0,1,\ldots,M;\quad t_k:=k\tau,\quad k=0,1,2,\ldots.$$
Denote
\begin{align*}
&X_M:=\left\{v=\left(v_0,v_1,\ldots,v_M\right)^T\  | \ v_0=v_M\right\} \subseteq \mathbb{C}^{M+1},\\
&Y_M:=\mathrm{span}\left\{e^{i\mu_l(x-a)},\quad l=-M/2,\cdots, M/2-1\right\}.
\end{align*}
For any function $\psi(x)\in L^2(\Omega)$ and $\phi(x)\in C_0(\overline\Omega)$ or vector $\phi=(\phi_0,\phi_1,\cdots,\phi_M)^T\in X_M$, let $\mathcal{P}_M: L^2(\Omega)\rightarrow Y_M$ be the standard $L^2$-projection operator, and $I_M: C_0(\overline\Omega)\rightarrow Y_M$ or $I_M: X_M\rightarrow Y_M$ be the standard interpolation operator as
\be\label{opd}
(P_M \psi)(x)=\sum \limits_{l=-M/2}^{M/2-1} \wh{\psi}_le^{i\mu_l(x-a)}, \quad (I_M \phi)(x)=\sum \limits_{l=-M/2}^{M/2-1} \widetilde{\phi}_le^{i\mu_l(x-a)},
\ee
where $\wh{\psi}_l$ and $\wt{\phi}_l$ are the Fourier and discrete Fourier transform coefficients of the function $\psi(x)$ and vector $\phi$ (with $\phi_j=\phi(x_j)$ when involved), respectively, defined as
\be\label{trans}
\wh{\psi}_l=\fl{1}{b-a}\int_a^b \psi(x)e^{-i\mu_l(x-a)}dx, \quad
\wt{\phi}_l=\fl{1}{M}\sum\limits_{j=0}^{M-1}\phi_j e^{-i\mu_l(x_j-a)},\quad l=-\fl{M}{2},\cdots, \fl{M}{2}-1.
\ee
Concerning the projection and interpolation operators, we review the standard estimates for the errors.
\begin{lemma}\cite{shen} \label{lemma}
For any $0 \leq \mu \leq k$, we have
\be\label{proc-e}
\|v-P_M(v)\|_\mu \leq Ch^{k-\mu}\|v\|_{k},\quad  \|P_N(v)\|_k\le C\|v\|_k,\quad \forall v\in H_{\rm p}^k(\Omega).
\ee
Moreover, if $k>1/2$, we have
\be\label{inter-e}
\|v-I_M(v)\|_{\mu} \leq C h^{k-\mu}\|v\|_{k},\quad  \|I_M(v)\|_k\le C\|v\|_k,\quad  \forall v\in H_{\rm p}^k(\Omega).
\ee
Here $C>0$ is a generic constant independent of $h$ and $v$.
\end{lemma}

A Fourier pseudospectral method for discretizing \eqref{gB} is to find
\be\label{sp}
z_{_M}(x,t)=\sum\limits_{l=-M/2}^{M/2-1} \wh{z}_l(t)e^{i\mu_l(x-a)},
\ee
such that
\be
\partial_{tt} z_{_M}-\partial_x^2 z_{_M}+\partial_x^4 z_{_M}-\partial_{xx} P_M(f(z_{_M}))=0,  \quad x \in \Og,\quad t>0. \label{GBN}
 \ee
Substituting \eqref{sp} into \eqref{GBN} and noticing the orthogonality of $\{e^{i\mu_l(x-a)}: -M/2\le l<M/2\}$, we obtain around time $t_k=k\tau$ $(k\ge0)$
\be
 \fl{d^2}{ds^2} \wh{z}_l(t_k+s)+\tht_l^2\wh{z}_l(t_k+s)+\mu_l^2\wh{\rho}_l(t_k+s)=0, \label{ZakNl}
 \ee
where $\rho=f(z_{_M})$ and $\tht_l=\sqrt{\mu_l^2+\mu_l^4}$. Applying the variation-of-constants formula \cite{zhao2016error} for $k\ge 0$ and $s\in \mathbb{R}$, the general solution of \eqref{ZakNl} can be written as follows for any $s\in\mathbb{R}$,
\be\label{forsol}
\begin{split}
&\wh{z}_0(t_k+s)=\wh{z}_0(t_k)+s\wh{z}'_0(t_k),\\
&\wh{z}_l(t_k+s)=\cos {(\tht_l s)}\wh{z}_l(t_k)+\fl{\sin (\tht_l s)}{\tht_l}\wh{z}'_l(t_k)-
\fl{\mu_l^2}{\tht_l}\int_0^s\wh{\rho}_l(t_k+\og)\sin{(\tht_l(s-\og))}d\og,\,\,\,\, l\neq 0.
\end{split}
\ee
Differentiating \eqref{forsol} with respect to $s$, we obtain
\be\label{forsol'}
\begin{split}
&\wh{z}_0'(t_k+s)=\wh{z}_0'(t_k),\\
&\wh{z}_l'(t_k+s)=-\tht_l \sin {(\tht_l s)}\wh{z}_l(t_k)+\cos (\tht_l s)\wh{z}'_l(t_k)-\mu_l^2\int_0^s\wh{\rho}_l(t_k+\og)\cos{(\tht_l(s-\og))}d\og,\quad l\neq 0.
\end{split}
\ee
Evaluating \eqref{forsol} and \eqref{forsol'} with $s=\tau$ and approximating the integral by the trapezoid rule or the Deuflhard-type quadrature \cite{Deuf,zhao2016error}, we immediately get
\begin{align*}
&\wh{z}_0(t_{k+1})=\wh{z}_0(t_k)+\tau\wh{z}'_0(t_k),\quad \wh{z}_0'(t_{k+1})=\wh{z}_0'(t_k),\\
&\wh{z}_l(t_{k+1})\approx \cos {(\tht_l \tau)}\wh{z}_l(t_k)+\fl{\sin (\tht_l \tau)}{\tht_l}\wh{z}'_l(t_k)-
\fl{\tau \mu_l^2}{2\tht_l}\wh{\rho}_l(t_k),\quad l\neq 0,\\
&\wh{z}_l'(t_{k+1})\approx -\tht_l \sin {(\tht_l \tau)}\wh{z}_l(t_k)+\cos (\tht_l \tau)\wh{z}'_l(t_k)-\fl{\tau\mu_l^2}{2}\left[
\cos(\tht_l\tau)\wh{\rho}_l(t_{k})+
\wh{\rho}_l(t_{k+1})\right],\quad l\neq 0.
\end{align*}

For implementation, the integrals for computing the Fourier transform coefficients are usually approximated by the numerical quadratures of \eqref{trans}. Let $z_j^k$ and $\dot{z}_j^k$ be the approximations of $z(x_j,t_k)$ and $\partial_t z(x_j,t_k)$, respectively, for
$0\le j<M$ and $k\ge0$; and denote $z^k$ and $\dot{z}^k$ as the vectors with components $z_j^k$ and $\dot{z}_j^k$, respectively. Choosing $z_j^0=z_0(x_j)$, $\dot{z}_j^0=z_1(x_j)$ for $0\le j<M$, a Fourier pseudospectral discretizatin for the problem \eqref{gB} reads
\be\label{sch1}
z_j^{k+1}=\sum\limits_{l=-M/2}^{M/2-1} \wt{z^{k+1}_l}e^{i\mu_l(x_j-a)},\quad
\dot{z}_j^{k+1}=\sum\limits_{l=-M/2}^{M/2-1} \wt{\dot{z}^{k+1}_l}e^{i\mu_l(x_j-a)},
\ee
where
\be\label{sch}
\begin{split}
&\wt{z^{k+1}_0}=\wt{z^{k}_0}+\tau \wt{\dot{z}^k_0},\quad \wt{\dot{z}^{k+1}_0}=\wt{\dot{z}^{k}_0},\\
&\wt{z^{k+1}_l}=\cos {(\tht_l \tau)}\wt{z^k_l}+\fl{\sin (\tht_l \tau)}{\tht_l}\wt{\dot{z}^k_l}-
\fl{\tau \mu_l^2}{2\tht_l}\sin {(\tht_l \tau)}\wt{(f(z^k))}_l,\quad l\neq 0,\\
&\wt{\dot{z}^{k+1}_l}=-\tht_l \sin {(\tht_l \tau)}\wt{z^k_l}+\cos (\tht_l \tau)\wt{\dot{z}^k_l}-\fl{\tau\mu_l^2}{2}\left[
\cos(\tht_l\tau)\wt{(f(z^k))}_l+
\wt{(f(z^{k+1}))}_l\right], \quad l\neq 0.
\end{split}
\ee
The above scheme is clearly explicit and very efficient due to the fast  discrete Fourier transform. The memory cost is $O(M)$ and the computational cost per time step is $O(M\ln(M))$.

\section{Convergence analysis}
For simplicity of notation, we denote the trigonometric interpolations of numerical solutions of \eqref{sch1} as
\[z_I^k(x):=I_M(z^k)(x),\quad \dot{z}_I^k(x):=I_M(\dot{z}^k)(x), \quad x\in \Omega.\]
Define the error functions as
\[e^k(x):=z(x,t_k)-z_I^k(x),\quad \dot{e}^k(x):=\partial_t z(x,t_k)-\dot{z}_I^k(x),\quad x\in\Omega,\quad k=0, 1, \ldots.\]
Then we have the following error estimates for \eqref{sch1} with \eqref{sch}.
\begin{theorem}\label{main}
Suppose $m>1/2$ and $\sigma\ge 4$. Let the solution of the GB equation \eqref{gB}  satisfies the regularity
properties $z\in C([0,T]; H^{m+\sigma}_{\mathrm{p}})\bigcap C^1([0,T]; H^{m+\sigma-2}_{\mathrm{p}})\cap
C^2([0,T]; H^{m}_{\mathrm{p}})$.  There exist $\tau_0>0$, $h_0>0$ sufficiently small such that when $\tau\le \tau_0$ and $h\le h_0$, we have the following error estimate for the numerical scheme \eqref{sch1} with \eqref{sch}:
\be\label{eror}
\|{e}^n\|_{m}+\|\dot{e}^n\|_{m-2}\leq K (\tau^2+h^{\sigma}),  \quad  0\le n\le T/\tau.
\ee
Furthermore, we have
\be\label{bound}
\|z_I^n\|_{m}\le K_1+1,\quad \|\dot{z}_I^n\|_{m-2}\le K_2+1,
\ee
where $K_1:=\|z\|_{L^\infty([0,T]; H^m)}$ and $K_2:=\|\p_t z\|_{L^\infty([0,T]; H^{m-2})}$.
\end{theorem}

Before giving the proof of Theorem \ref{main}, we present some auxiliary results established in \cite{Super}, which will be applied repeatedly in our proof.
\begin{proposition}\cite{Super}\label{Super}
For any function $g\in C^\infty(\mathbb{C}, \mathbb{C})$ and $s >1/2$, there exists a nondecreasing function $\chi_g: \mathbb{R}^+\rightarrow \mathbb{R}^+$ such that
\be\label{ct}
\|g(u)\|_s \le \|g(0)\|_s+\chi_g(\|u\|_{L^\infty})\|u\|_s,\quad \forall u\in H^s.
\ee
For all $v, w\in B_R^s:=\{u\in H^s: \|u\|_s\le R\}$, we have
\be\label{Lip}
\|g(v)-g(w)\|_s\le \alpha(g, R)\|v-w\|_s,
\ee
where $\alpha(g, R)=\|g'(0)\|_s+R\chi_{g'}(cR)$ is nondecreasing with respect to $R$, with $c>0$ being the constant for the Sobolev imbedding $\|\cdot\|_{L^\infty}\le c\|\cdot\|_s$.
\end{proposition}

\emph{Proof of Theorem \ref{main}.} We deduce \eqref{eror} and \eqref{bound} by induction. For $n=0$, \eqref{eror} is obvious by noticing
\[\|e^0\|_m=\|z_0-I_M(z_0)\|_m\le Ch^\sigma \|z_0\|_{m+\sigma},\quad \|\dot{e}^0\|_{m-2}=\|z_1-I_M(z_1)\|_{m-2}\le Ch^\sigma \|z_1\|_{m+\sigma-2}.\]
Furthermore, \eqref{bound} can be obtained by the triangle inequality when $\tau$ and $h$ are small enough.

Now assume \eqref{eror} and \eqref{bound} are valid for $n=0, \ldots, k<T/\tau$, next we show \eqref{eror} and \eqref{bound} are true for $k+1$. To precede, denote the projected error functions as
\[e_M^n(x):=P_M(e^n)(x)=\sum\limits_{l=-M/2}^{M/2-1} \wh{e^n_l}e^{i\mu_l(x-a)},\quad \dot{e}_M^n(x):=P_M(\dot{e}^n)(x)=\sum\limits_{l=-M/2}^{M/2-1}\wh{\dot{e}^n_l}e^{i\mu_l(x-a)},\]
where the corresponding coefficients satisfy
\[\wh{e^n_l}=\wh{z}_l(t_n)-\wt{z^n_l},\quad
\wh{\dot{e}^n_l}=\wh{ z}'_l(t_n)-\wt{\dot{z}^n_l},\quad l=-M/2,\ldots,M/2-1.\]
Applying the triangle inequality and Lemma \ref{lemma}, we have
\begin{align*}
&\|e^n\|_m\le \|e_M^n\|_m+\|P_M(z(\cdot, t_n))-z(\cdot, t_n)\|_m\le  \|e_M^n\|_m+Ch^\sigma \|z(\cdot, t_n)\|_{m+\sigma},\\
&\|\dot{e}^n\|_{m-2}\le \|\dot{e}_M^n\|_{m-2}+\|P_M(\p_tz(\cdot, t_n))-\p_tz(\cdot, t_n)\|_{m-2}\le  \|\dot{e}_M^n\|_{m-2}+Ch^\sigma \|\p_t z(\cdot, t_n)\|_{m+\sigma-2},
\end{align*}
thus it suffices to show \eqref{bound} and
\be\label{emn}
\|e_M^n\|_m+\|\dot{e}_M^n\|_{m-2}\le K (\tau^2+h^{\sigma}) .
\ee
Define the local truncation errors as
\[\xi^n(x):=\sum\limits_{l=-M/2}^{M/2-1}\wh{\xi^n_l}e^{i\mu_l(x-a)},\quad
\dot{\xi}^n(x):=\sum\limits_{l=-M/2}^{M/2-1}
\wh{\dot{\xi}^n_l}e^{i\mu_l(x-a)},\]
where $\wh{\xi^n_0}=\wh{\dot{\xi}^n_0}=0$, and for $l\neq 0$,
\begin{align*}
\wh{\xi^n_l}&=\wh{z}_l(t_{n+1})-\wh{z}_l(t_n)\cos(\tht_l\tau)-
\fl{\sin(\tht_l\tau)}{\tht_l}\wh{z}_l'(t_n)+\fl{\tau\mu_l^2}{2\tht_l}
\sin(\tht_l\tau)\wh{\rho}_l(t_n),\\
\wh{\dot{\xi}^n_l}&=\wh{z}'_l(t_{n+1})+\tht_l\sin(\tht_l\tau)\wh{z}_l(t_n)-
\cos(\tht_l\tau)\wh{z}_l'(t_n)+\fl{\tau\mu_l^2}{2}\left[\cos(\tht_l\tau)\wh{\rho}_l(t_n)+\wh{\rho}_l(t_{n+1})\right],
\end{align*}
with $\rho(x,s)=f(z(x,s))$.
Adding the local truncation errors from the scheme, we are led to the error equations as
\be\label{ek}
\begin{split}
&\wh{e^{n+1}_0}=\wh{e^n_0}+\tau\wh{\dot{e}^n_0},\quad \wh{\dot{e}^{n+1}_0}=\wh{\dot{e}^{n}_0},\\
&\wh{e^{n+1}_l}=\cos(\tht_l\tau)\wh{e^n_l}+\fl{\sin(\tht_l\tau)}{\tht_l}\wh{\dot{e}^n_l}+
\wh{\xi^n_l}-\fl{\tau\mu_l^2}{2\tht_l}\sin(\tau\tht_l)\wh{\eta^n_l},\quad l\neq 0,\\
&\wh{\dot{e}^{n+1}_l}=-\tht_l\sin(\tht_l\tau)\wh{e^n_l}+
\cos(\tht_l\tau)\wh{\dot{e}^n_l}+\wh{\dot{\xi}^n_l}-\fl{\tau\mu_l^2}{2}
\left[\cos(\tht_l\tau)\wh{\eta^n_l}+\wh{\eta^{n+1}_l}\right], \quad l\neq 0,
\end{split}
\ee
where
\[\wh{\eta^n_l}=\wh{\rho}_l(t_n)-\wt{\rho^n_l} ,\quad l\neq 0,\quad \rho^n=f(z^n).\]
It follows from \eqref{ek} that
\begin{align*}
|\wh{e^{n+1}_l}|^2&\le (1+\tau)\left|\cos(\tht_l\tau)\wh{e^n_l}+\fl{\sin(\tht_l\tau)}{\tht_l}\wh{\dot{e}^n_l}\right|^2+(1+\fl{1}{\tau})
\left|\wh{\xi^n_l}+\fl{\tau\mu_l^2}{2\tht_l}\sin(\tau\tht_l)\wh{\eta^n_l}\right|^2,\\
|\wh{\dot{e}^{n+1}_l}|^2&\le(1+\tau)\left|\tht_l\sin(\tht_l\tau)\wh{e^n_l}-
\cos(\tht_l\tau)\wh{\dot{e}^n_l}\right|^2\\
&\quad+(1+\fl{1}{\tau})\left|\wh{\dot{\xi}^n_l}+\fl{\tau\mu_l^2}{2}
\left[\cos(\tht_l\tau)\wh{\eta^n_l}+\wh{\eta^{n+1}_l}\right]\right|^2,\quad l\neq 0,
\end{align*}
which yields for $l\neq 0$,
\begin{align}
\tht_l^2|\wh{e^{n+1}_l}|^2+|\wh{\dot{e}^{n+1}_l}|^2&\le (1+\tau)\left[\tht_l^2|\wh{e^{n}_l}|^2+|\wh{\dot{e}^{n}_l}|^2\right]+
(1+\fl{1}{\tau})\left[\tht_l^2\left|\wh{\xi^n_l}+\fl{\tau\mu_l^2}{2\tht_l}\sin(\tau\tht_l)
\wh{\eta^n_l}\right|^2\right.\nn\\
&\qquad\left.+\left|\wh{\dot{\xi}^n_l}+\fl{\tau\mu_l^2}{2}
\left[\cos(\tht_l\tau)\wh{\eta^n_l}+\wh{\eta^{n+1}_l}\right]\right|^2\right]\nn\\
&\le (1+\tau)\left[\tht_l^2|\wh{e^{n}_l}|^2+|\wh{\dot{e}^{n}_l}|^2\right]+
(1+\fl{1}{\tau})\Big[2\tht_l^2|\wh{\xi^n_l}|^2+\tau^2\mu_l^4
|\wh{\eta^n_l}|^2\Big.\nn\\
&\qquad\Big.+2|\wh{\dot{\xi}^n_l}|^2+\tau^2\mu_l^4\Big(|\wh{\eta^n_l}|^2+|\wh{\eta^{n+1}_l}|^2\Big)\Big].\label{enp}
\end{align}
Denote
\[\mathcal{E}^n=|\wh{e^n_0}|^2+|\wh{\dot{e}^n_0}|^2+\sum\limits_{-\fl{M}{2}\le l<\fl{M}{2},\,\, l\neq 0}(1+|\mu_l|^2)^{m-2}(\tht_l^2
|\wh{e^n_l}|^2+|\wh{\dot{e}^n_l}|^2).\]
By definition \eqref{sn}, we have
\be\label{Ens}
\mathcal{E}^n\sim \|e_M^n\|_m^2+\|\dot{e}_M^n\|_{m-2}^2,
\ee
where $p \lesssim q$ means there exist a constant $C>0$ such that $p\le C q$ and $p\sim q$ represents $q\lesssim p\lesssim q$.
Noticing that
\[|\wh{e^{n+1}_0}|^2\le (1+\tau)|\wh{e^n_0}|^2+(1+\fl{1}{\tau})\tau^2|\wh{\dot{e}^n_0}|^2,\quad |\wh{\dot{e}^{n+1}_0}|^2=|\wh{\dot{e}^{n}_0}|^2,\]
this together with \eqref{enp} derives that
\be\label{induc}
\mathcal{E}^{n+1}-\mathcal{E}^n \lesssim \tau \mathcal{E}^n+\fl{1}{\tau}\left(\|\xi^n\|_m^2+\|\dot{\xi}^n\|_{m-2}^2\right)
+\tau\left(\|\eta^n\|_m^2+\|\eta^{n+1}\|_m^2\right).
\ee
Next we estimate the local error and the error for the nonlinear term, respectively.

\emph{Local error.}
Applying the quadrature error
\[\int_0^\tau g(s)ds-\fl{\tau}{2}(g(0)+g(\tau))=-\fl{1}{2}\int_0^\tau g''(s)s(\tau-s)ds,\]
we get for $l\neq 0$,
\begin{align*}
\wh{\xi^n_l}&=\fl{\mu_l^2}{2\tht_l}\int_0^\tau s(\tau-s)[\sin(\tht_l(\tau-s))\wh{\rho}_l(t_k+s)]''ds
=\fl{\mu_l^2}{2\tht_l}\int_0^\tau s(\tau-s)P_l^n(s)ds,\\
\wh{\dot{\xi}^n_l}&=\fl{\mu_l^2}{2}\int_0^\tau s(\tau-s)[\cos(\tht_l(\tau-s))\wh{\rho}_l(t_k+s)]''ds
=\fl{\mu_l^2}{2}\int_0^\tau s(\tau-s)Q_l^n(s)ds,
\end{align*}
where
\begin{align*}
P_l^n(s)&=-\tht_l^2\sin(\tht_l(\tau-s))\wh{\rho}_l(t_n+s)-
2\tht_l\cos(\tht_l(\tau-s))\wh{\rho}_l'(t_n+s)+\sin(\tht_l(\tau-s))
\wh{\rho}''_l(t_n+s),\\
Q_l^n(s)&=-\tht_l^2\cos(\tht_l(\tau-s))\wh{\rho}_l(t_n+s)-2\tht_l
\sin(\tht_l(\tau-s))\wh{\rho}'_l(t_n+s)+\cos(\tht_l(\tau-s))\wh{\rho}''_l(t_n+s).
\end{align*}
Applying the H\"older's inequality, we get for $l\neq 0$,
\begin{align*}
|\wh{\xi^n_l}|&\le \fl{\mu_l^2}{2\tht_l}\left(\int_0^\tau s^2(\tau-s)^2ds\right)^{1/2}\left(\int_0^\tau |P_l^n(s)|^2ds\right)^{1/2}\le
\fl{\tau^{5/2}}{2}\left(\int_0^\tau |P_l^n(s)|^2ds\right)^{1/2},\\
|\wh{\dot{\xi}^n_l}|&\le \fl{\mu_l^2}{2}\left(\int_0^\tau s^2(\tau-s)^2ds\right)^{1/2}\left(\int_0^\tau |Q_l^n(s)|^2ds\right)^{1/2}\le
\fl{\mu_l^2\tau^{5/2}}{2}\left(\int_0^\tau |Q_l^n(s)|^2ds\right)^{1/2}.
\end{align*}
Thus
\begin{align*}
\|\xi^n\|^2_m&= \sum\limits_{l=-M/2}^{M/2-1}(1+|\mu_l|^2)^m|\wh{\xi^n_l}|^2\\
&\lesssim \tau^5 \sum\limits_{-\fl{M}{2}\le l<\fl{M}{2},\,\, l\neq 0}(1+|\mu_l|^2)^m
\int_0^\tau |P_l^n(s)|^2ds\\
&\lesssim \tau^5 \sum\limits_{-\fl{M}{2}\le l<\fl{M}{2},\,\, l\neq 0}(1+|\mu_l|^2)^m\left[\tht_l^4\int_0^\tau
|\wh{\rho}_l(t_n+s)|^2ds+\tht_l^2\int_0^\tau|\wh{\rho}'_l(t_n+s)|^2ds\right.\\
&\qquad\qquad\qquad\qquad\qquad\qquad\qquad\left.+\int_0^\tau |\wh{\rho}''_l(t_n+s)|^2ds\right]\\
&\lesssim \tau^{5}\int_0^\tau \left[\|\rho(\cdot,t_n+s)\|_{m+4}^2+\|\p_t\rho(\cdot, t_n+s)\|_{m+2}^2+
\|\p_t^2\rho(\cdot,t_n+s)\|_{m}^2\right]ds.
\end{align*}
Noticing that
\[\p_t \rho=f'(z)\p_tz,\quad \p_t^2\rho=f''(z)(\p_tz)^2+f'(z)\p_t^2z,\]
recalling $m>1/2$, using the bilinear estimate \cite{Adams}
\be\label{bi}
\|f g\|_{r}\le C_r \|f\|_r\|g\|_r,\quad r>1/2,
\ee
and Proposition \ref{Super}, we yield that
\begin{align*}
&\|\rho(\cdot,t)\|_{m+4}\le \|f(0)\|_{m+4}+\chi_f(\|z(\cdot,t)\|_{L^\infty})\|z(\cdot, t)\|_{m+4}\le
\|f(0)\|_{m+4}+\chi_f(L)\|z(\cdot, t)\|_{m+4},\\
&\|\p_t\rho(\cdot, t)\|_{m+2}\le C_m(\|f'(0)\|_{m+2}+\chi_{f'}(L)\|z(\cdot, t)\|_{m+2})\|\p_tz(\cdot,t)\|_{m+2},\\
&\|\p_t^2\rho(\cdot, t)\|_{m}\le C_m(\|f'(0)\|_{m}+\chi_{f'}(L)\|z(\cdot, t)\|_{m})\|\p_t^2z(\cdot,t)\|_{m}\\
&\qquad\qquad\qquad+ C_m^2 (\|f''(0)\|_{m}+\chi_{f''}(L)\|z(\cdot,t)\|_{m})\|\p_tz(\cdot,t)\|^2_{m},
\end{align*}
where $L=\|z\|_{L^\infty([0,T]; L^\infty)}$, which is finite by recalling that $z\in C([0,T]; H^{m+\sigma}_{\mathrm{p}})$. By the assumptions on the regularities of $z$, we get that
\[\|\xi^n\|_m^2\lesssim \tau^6.\]
Similar approach leads to that
\begin{align*}
\|\dot{\xi}^n\|_{m-2}^2&= \sum\limits_{l=-M/2}^{M/2-1}(1+|\mu_l|^2)^{m-2}|\wh{\dot{\xi}^n_l}|^2\\
&\lesssim \tau^5 \sum\limits_{-\fl{M}{2}\le l<\fl{M}{2},\,\, l\neq 0}\mu_l^4(1+|\mu_l|^2)^{m-2}
\int_0^\tau |Q_l^n(s)|^2ds\\
&\lesssim \tau^5 \sum\limits_{-\fl{M}{2}\le l<\fl{M}{2},\,\, l\neq 0}(1+|\mu_l|^2)^{m}\left[\tht_l^4\int_0^\tau
|\wh{\rho}_l(t_n+s)|^2ds+\tht_l^2\int_0^\tau|\wh{\rho}'_l(t_n+s)|^2ds\right.\\
&\quad\left.
+\int_0^\tau |\wh{\rho}''_l(t_n+s)|^2ds\right]\lesssim \tau^{6}.
\end{align*}
Hence we get the estimate on the local errors
\be\label{localer}
\|\xi^n\|^2_m+\|\dot{\xi}^n\|_{m-2}^2\le M_1\tau^6,
\ee
where $M_1$ depends on $m$, $f$, $\|z\|_{L^\infty([0,T]; H^{m+4})}$,  $\|\p_tz\|_{L^\infty([0,T]; H^{m+2})}$ and $\|\p_t^2z\|_{L^\infty([0,T]; H^{m})}$.
\medskip

\emph{Error of nonlinear terms.}
By definition, we have
\begin{align}
\|\eta^n\|_m&=\|I_M(\rho^n)-P_M(\rho(\cdot,t_n))\|_{m}\nn\\
&\le \|I_M(f(z_I^n))-I_M(\rho(\cdot,t_n))\|_{m}+\|I_M(\rho(\cdot,t_n))-P_M(\rho(\cdot,t_n))\|_m\nn\\
&\lesssim \|f(z_I^n)-f(z(\cdot,t_n))\|_m+h^{\sigma}\|\rho(\cdot,t_n)\|_{m+\sigma}\nn\\
&\le \alpha(f, \max\{\|z_I^n\|_m, \|z(\cdot,t_n)\|_m\})\|z_I^n-z(\cdot,t_n)\|_m+h^{\sigma}
(\|f(0)\|_{m+\sigma}+\chi_f(L)\|z(\cdot,t_n)\|_{m+\sigma})\nn\\
&\lesssim \alpha(f, \max\{\|z_I^n\|_m, \|z(\cdot,t_n)\|_m\})(\|e_M^n\|_m+h^\sigma \|z(\cdot,t_n)\|_{m+\sigma})+h^{\sigma}.\label{etak}
\end{align}
Similar approach yields that
\be\label{etakp}
\|\eta^{n+1}\|_m\lesssim \alpha(f, \max\{\|z_I^{n+1}\|_m, \|z(\cdot,t_{n+1})\|_m\})(\|e_M^{n+1}\|_m+h^\sigma \|z(\cdot,t_{n+1})\|_{m+\sigma})+h^{\sigma}.
\ee
To estimate this, we need an a prior bound for $\|z_I^{n+1}\|_m$.
It follows from the scheme that
\begin{align}
\|z_I^{n+1}\|_m&=\left(\sum\limits_{l=-M/2}^{M/2-1}
(1+|\mu_l|^2)^m|\wt{z^{n+1}_l}|^2\right)^{1/2}\nn\\
&\le \left|\wt{z^{n}_0}+\tau\wt{\dot{z}^n_0}\right|+\|z_I^n\|_m+\left(\sum\limits_{-M/2\le l<M/2, l\neq 0}
(1+|\mu_l|^2)^m\tht_l^{-2}|\wt{\dot{z}^n_l}|^2\right)^{1/2}+\fl{\tau}{2}\|f(z_I^n)\|_m\nn\\
&\le 2\|z_I^n\|_m+\tau\left|\wt{\dot{z}^n_0}\right|+c_1\|\dot{z}_I^n\|_{m-2}+\fl{\tau}{2}\left(\|f(0)\|_m+\chi_f(\|z_I^n\|_{L^\infty})
\|z_I^n\|_m\right)\nn\\
&\le 2(K_1+1)+c_1(K_2+1)+\|f(0)\|_m+\chi_f(c(K_1+1))(K_1+1)+\left|\wt{\dot{z}^0_0}\right|,\label{zI}
\end{align}
when $\tau\le 1$. Here we have used the property \eqref{bi}, $\wt{\dot{z}^n_0}=\wt{\dot{z}^{k-1}_0}=\ldots=\wt{\dot{z}^0_0}$, the induction and the inequality
\[\tht_l^{-2}\le c_1^2(1+|\mu_l|^2)^{-2},\quad l\neq 1, \]
with $c_1^2=1+\fl{(b-a)^2}{4\pi^2}$.
Hence combining \eqref{etak}, \eqref{etakp} and \eqref{zI}, we get
\be\label{non}
\|\eta^n\|_m^2+\|\eta^{n+1}\|_m^2\le C(f, K_1)\mathcal{E}^n+C(f, K_1, K_2)\mathcal{E}^{n+1}+C(f, K_1, K_2, K_\sigma)h^{2\sigma}.
\ee

\medskip
Combining \eqref{induc}, \eqref{localer} and \eqref{non},
we derive that
\[\mathcal{E}^{n+1}-\mathcal{E}^n \lesssim  \tau^5+\tau h^{2\sigma}+\tau\left(\mathcal{E}^n+\mathcal{E}^{n+1}\right).\]
Summing the above inequality for $n=0, 1, \ldots, k$, one gets
\[\mathcal{E}^{k+1}-\mathcal{E}^{0}\lesssim \tau \sum\limits_{n=0}^{k+1} \mathcal{E}^n+\tau^4+h^{2\sigma}.\]
Applying the discrete Gronwall's inequality, when $\tau$ is sufficiently small, we have
\[\mathcal{E}^{k+1}\lesssim \tau^4+h^{2\sigma},\]
which immediately gives \eqref{emn} by recalling \eqref{Ens}. Finally \eqref{bound} can be obtained by \eqref{eror} and the triangle inequality.
\hfill $\square$ \bigskip

\section{Numerical experiments}
In this section, we first test the order of accuracy of the DEI-FP scheme \eqref{sch1} with \eqref{sch}. Then we apply this method to investigate some long time dynamics of the GB euqation. For all the numerical experiments, we choose the commonly used nonlinearity $f(z)=z^2$.
\subsection{Accuracy test}
In the first experiment, we test the convergence of the Deuflhard-type exponential integrator Fourier pseudospectral scheme \eqref{sch1} with \eqref{sch} for the solitary wave solution.

\medskip
{\sl Example 1.} The well-known soliton solution of the GB equation \eqref{GB} is given by \cite{Alex19,zhang2018}
\be\label{sol-ex}
z(x,t)=-A\,\sech^2\left(\sqrt{A/6}(x-vt-x_0)\right),\quad v=\pm \sqrt{1-2A/3},
\ee
where $A$, $x_0$ and $v$ represent the amplitude, the initial location and the velocity of the soliton, respectively.

Noticing that the solitary wave decays exponentially in the far field, this enables us to consider the GB equation on a bounded interval $[-a, a]$ with periodic boundary conditions when $a$ is large enough such that the artificial boundaries are located far out enough for the theoretical solution to satisfy the periodic boundary conditions except for a negligible remainder. Here we choose $A=3/8$, $x_0=0$ and the torus $\Omega=(-60,60)$. Denote $z^{\tau,h}$ and $\dot{z}^{\tau,h}$ as the numerical solutions obtained by the DEI-FP method with mesh size $h$ and time step $\tau$ for approximating the exact solutions $z(\cdot,t)$ and $\p_tz(\cdot,t)$. To quantify the numerical error, we define the error function as
\[e_m^{\tau,h}:=\|I_M(z^{\tau,h})-z(\cdot,t)\|_m+\|I_M(\dot{z}^{\tau,h})-\p_t z(\cdot,t)\|_{m-2}.\]

\begin{figure}[h!]
\begin{minipage}[t]{0.5\linewidth}
\centering
\includegraphics[height=4.9cm,width=7.1cm]{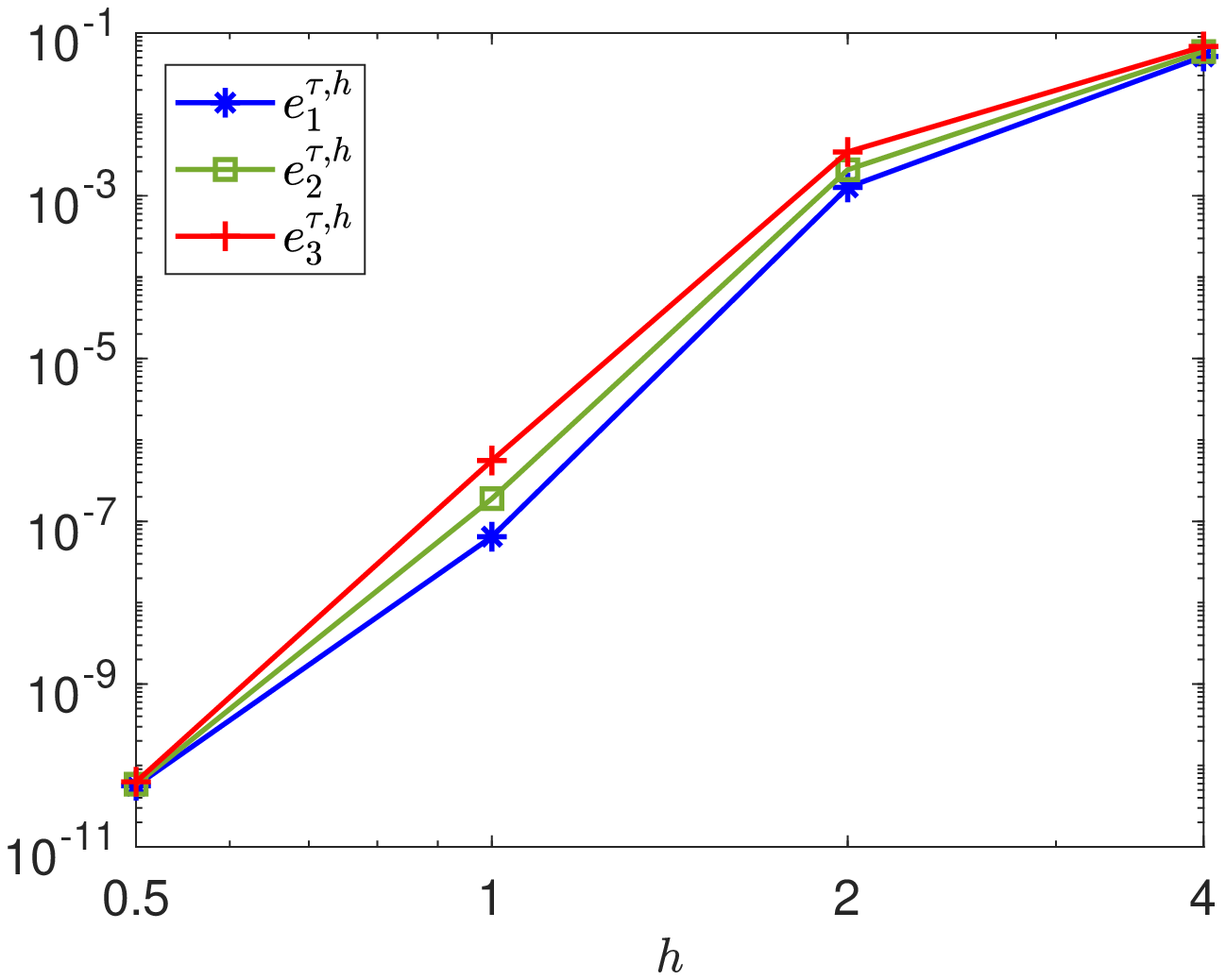}
\end{minipage}%
\hspace{3mm}
\begin{minipage}[t]{0.5\linewidth}
\centering
\includegraphics[height=4.9cm,width=7.1cm]{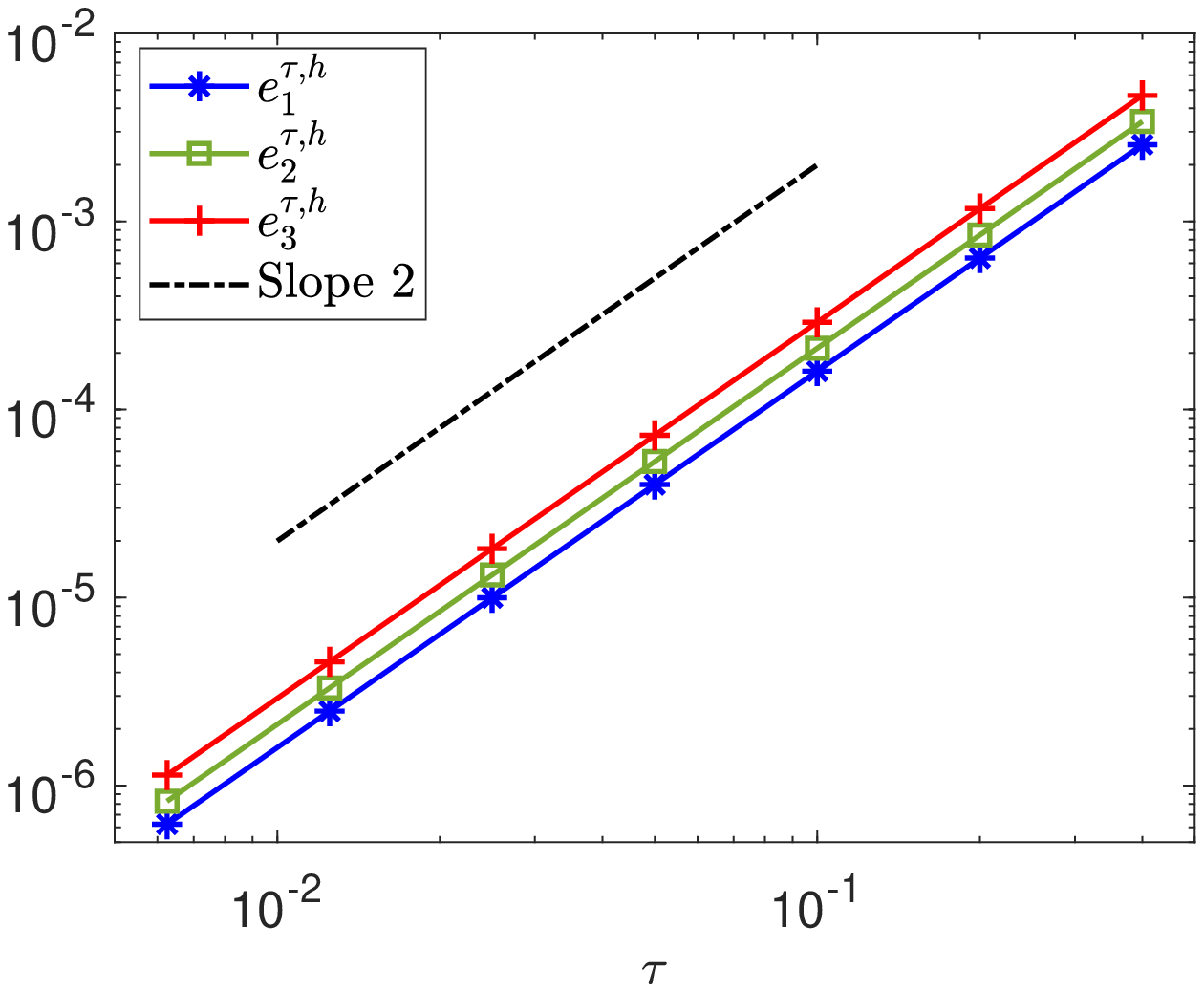}
\end{minipage}
\vspace{-3mm}
\caption{Spatial (left) and temporal (right) errors of the DEI-FP scheme for the soliton solution under different mesh size and time step size.}\label{sol}
\end{figure}

Fig. \ref{sol} displays the spatial and temporal errors of DEI-FP method \eqref{sch1}-\eqref{sch} for the solitary wave solution at $T=2$ under various values of $\tau$ and $h$. The errors are quantified in several norms with $m=1,2,3$. For spatial error analysis, we take a tiny time step $\tau=10^{-6}$ such that the temporal discretization error is negligible; for temporal error analysis, we set the mesh size $h=1/8$ such that the spatial error can be ignorable. It can be clearly observed that the scheme converges spectrally and quadratically in space and time, respectively, which confirms the theoretical result in Theorem \ref{main}.

\begin{figure}[h!]
\begin{minipage}[t]{0.5\linewidth}
\centering
\includegraphics[height=4.9cm,width=7.1cm]{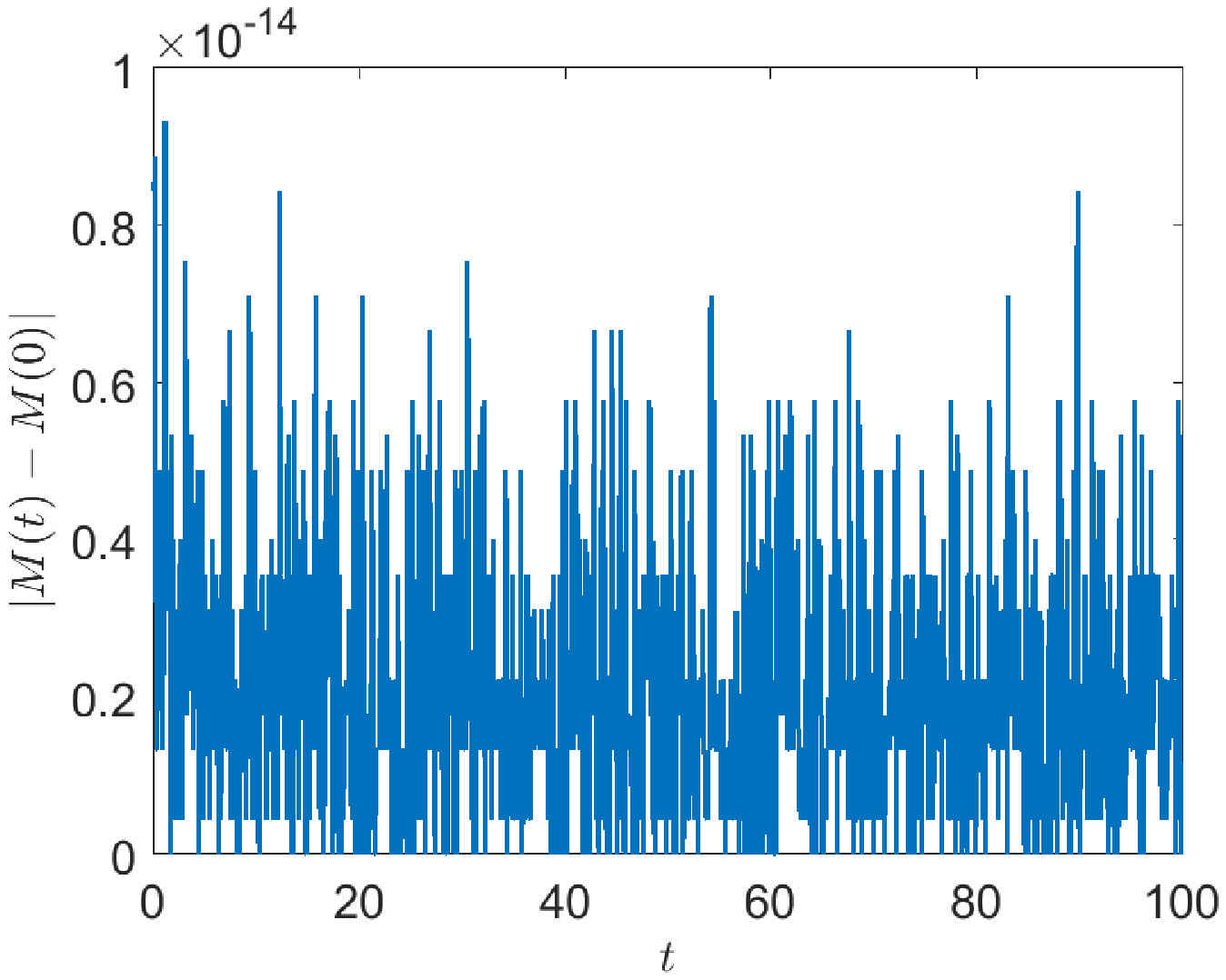}
\end{minipage}%
\hspace{3mm}
\begin{minipage}[t]{0.5\linewidth}
\centering
\includegraphics[height=4.9cm,width=7.1cm]{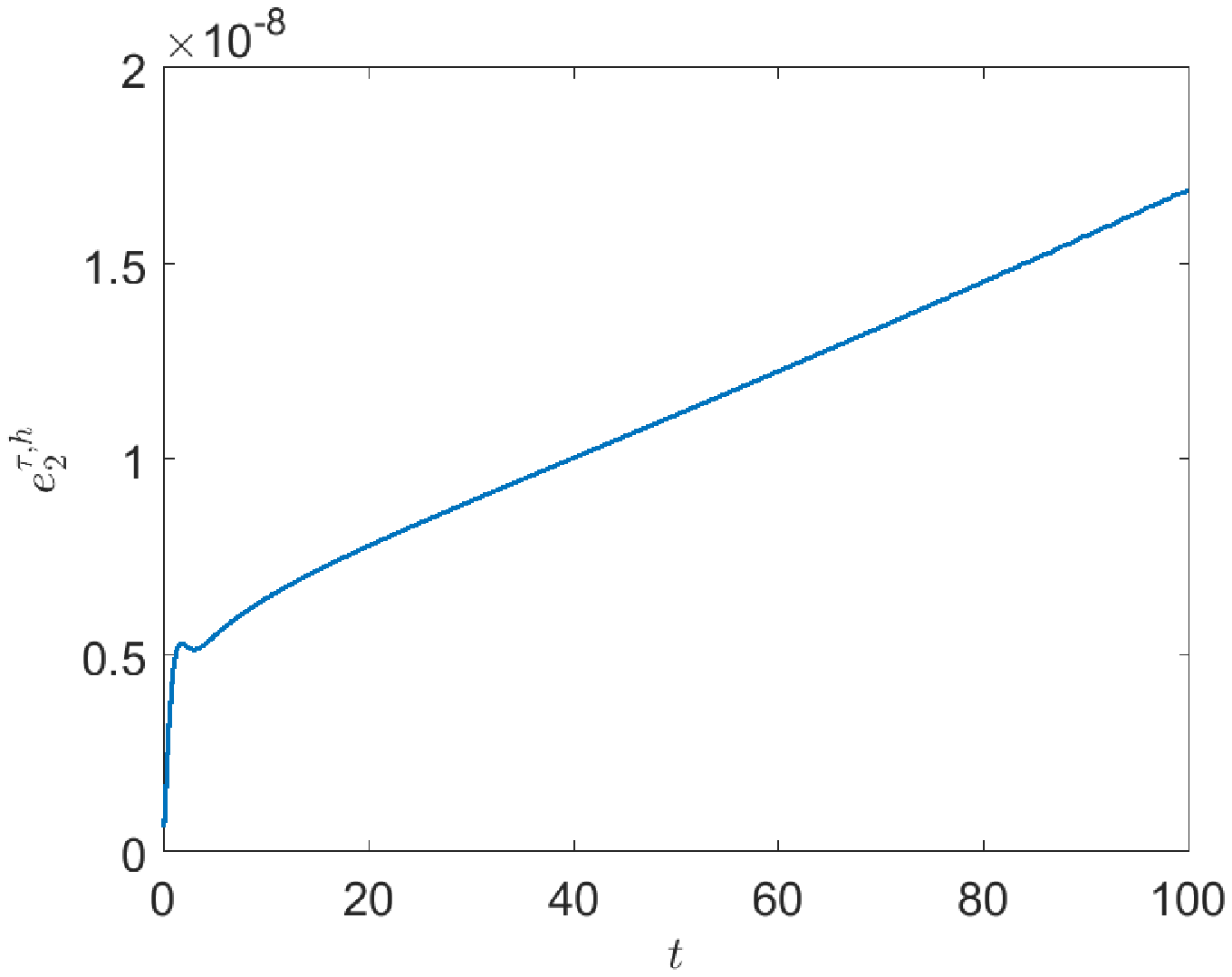}
\end{minipage}
\vspace{-3mm}
\caption{Conservation of mass (left) and long-time errors of the DEI-FP method (right).}\label{long}
\end{figure}

Next we investigate the long time behavior of the DEI-FP method. Noticing that for the soliton solution \eqref{sol-ex}, $\int_\Omega \p_t z(x, 0)dx=0$, this implies that the mass is conserved:
\[M(t):=\int_\Omega z(x,t)dx.\]
Fig. \ref{long} shows the conservation law of mass for the numerical solitary wave solution (left) and the long time error of the DEI-FP scheme (right). Here the solution is obtained on a bounded interval $\Omega=[-300, 300]$ under $h=1/8$ and $\tau=0.001$. We clearly see that the DEI-FP method is reliable and excellent for long time intervals.

\subsection{Birth of solitons}
{\sl Example 2.} In this experiment, we still use the initial data for the solitary wave solution \eqref{sol-ex}:
\be\label{ini_set2}
\begin{split}
&u_0(x)=-A\,\sech^2\left(\sqrt{A/6}\,x\right),\\
&u_1(x)=-A\sqrt{1-2A/3}\sqrt{2A/3}\,\sech^2(\sqrt{A/6}\,x)\tanh(\sqrt{A/6}\,x).
\end{split}
\ee

\begin{figure}[h!]
\begin{minipage}[t]{0.5\linewidth}
\centering
\includegraphics[height=3cm,width=7.1cm]{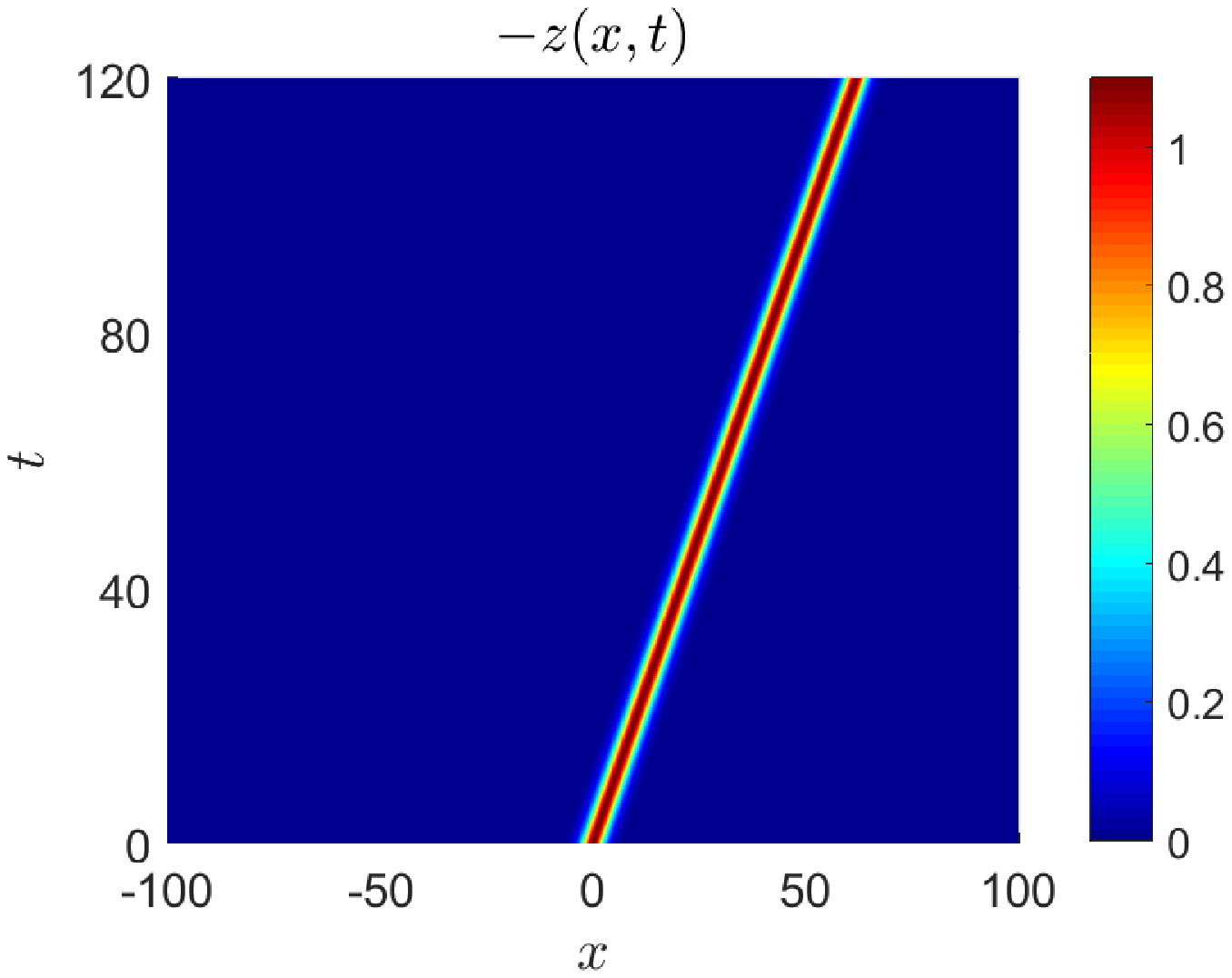}
\end{minipage}%
\hspace{3mm}
\begin{minipage}[t]{0.5\linewidth}
\centering
\includegraphics[height=3cm,width=7.1cm]{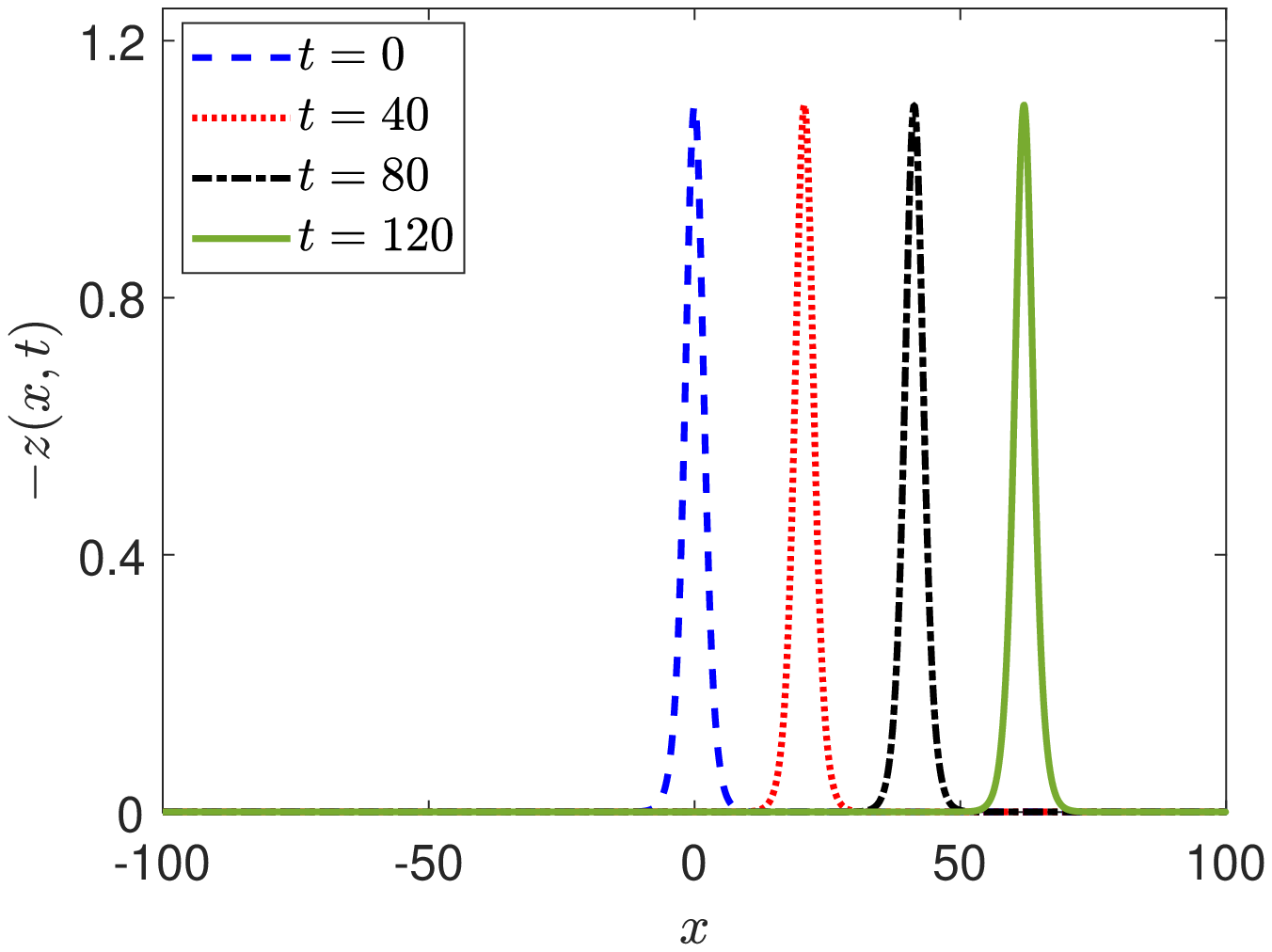}
\end{minipage}
\begin{minipage}[t]{0.5\linewidth}
\centering
\includegraphics[height=3cm,width=7.1cm]{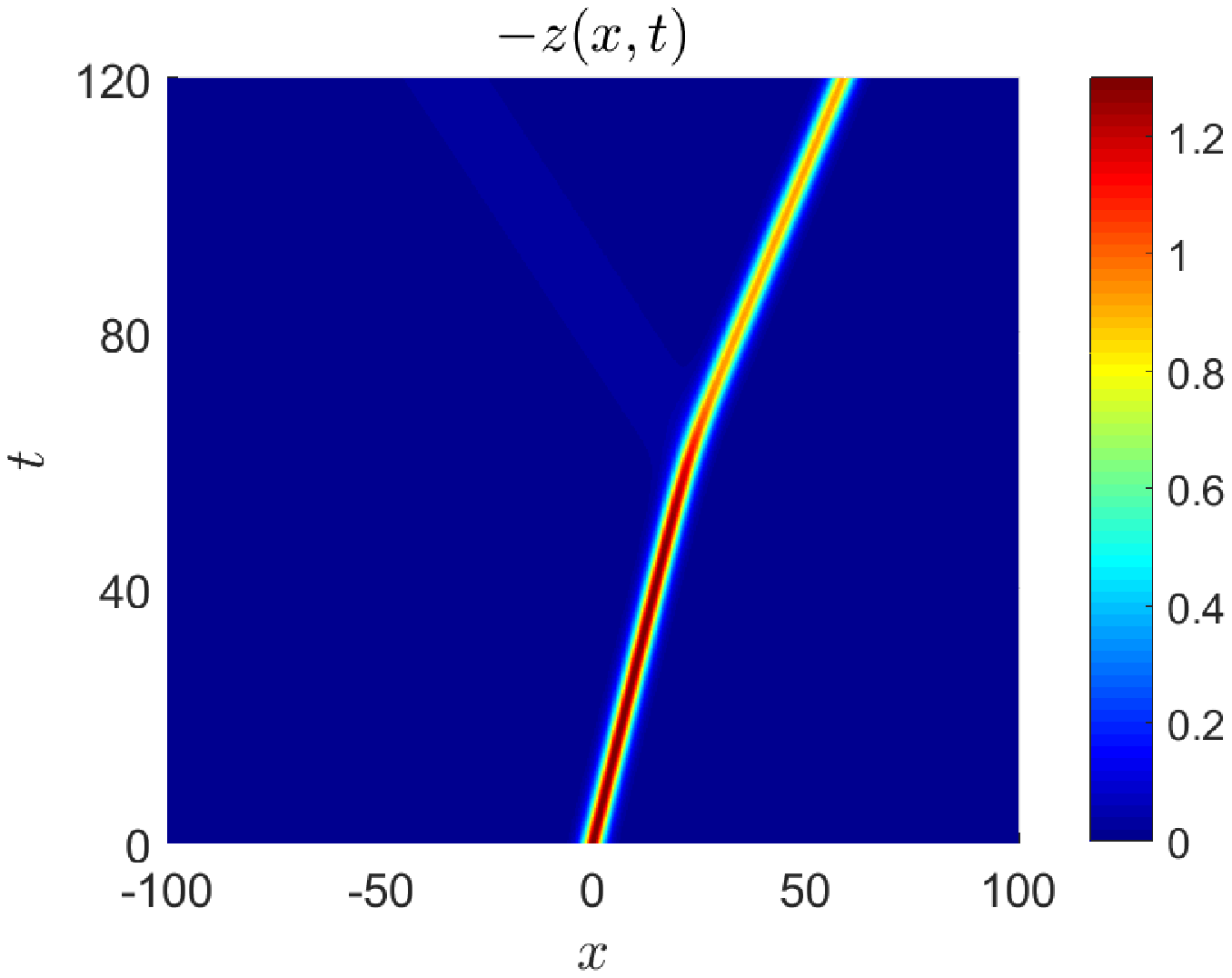}
\end{minipage}%
\hspace{3mm}
\begin{minipage}[t]{0.5\linewidth}
\centering
\includegraphics[height=3cm,width=7.1cm]{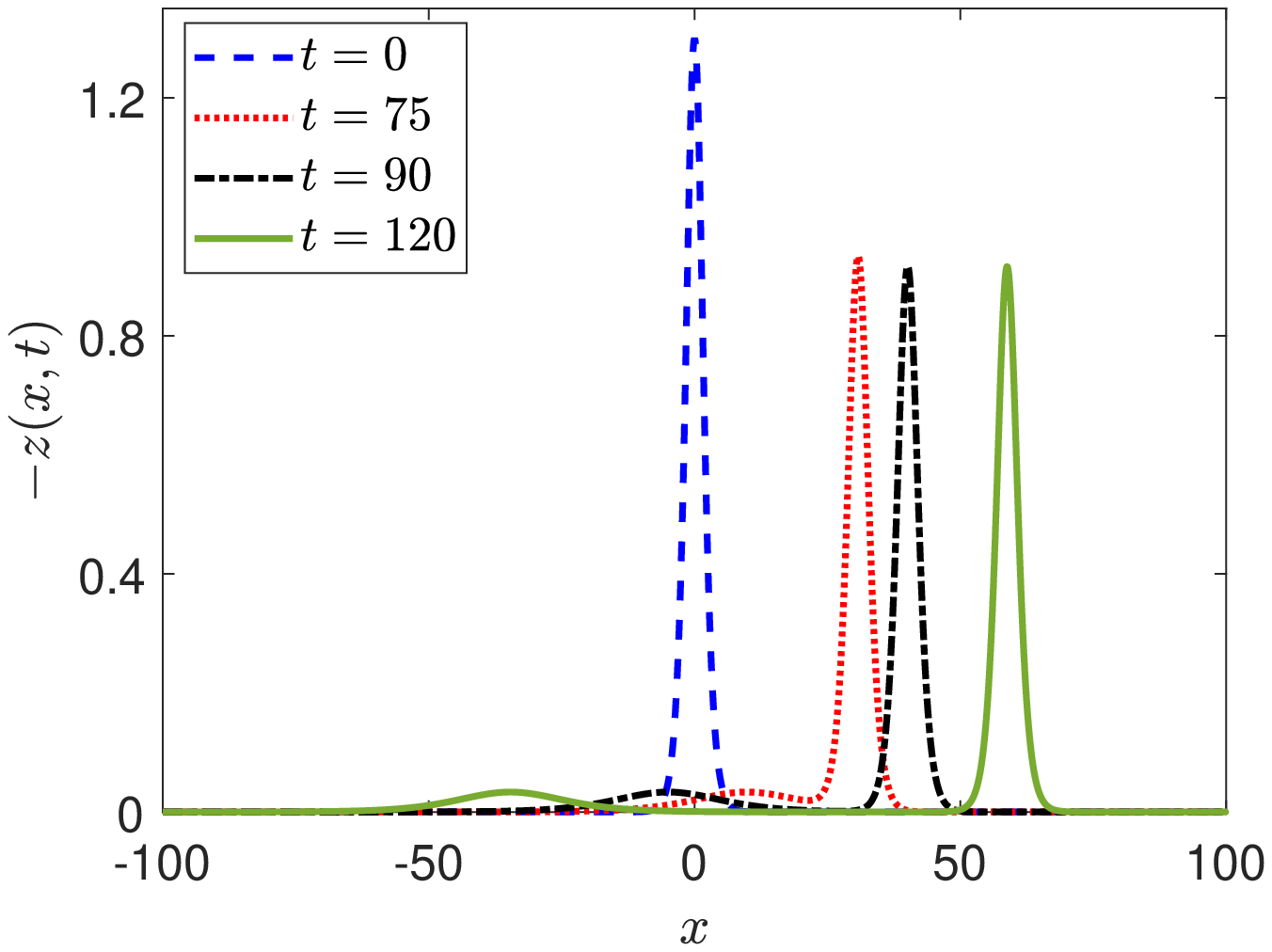}
\end{minipage}
\begin{minipage}[t]{0.5\linewidth}
\centering
\includegraphics[height=3cm,width=7.1cm]{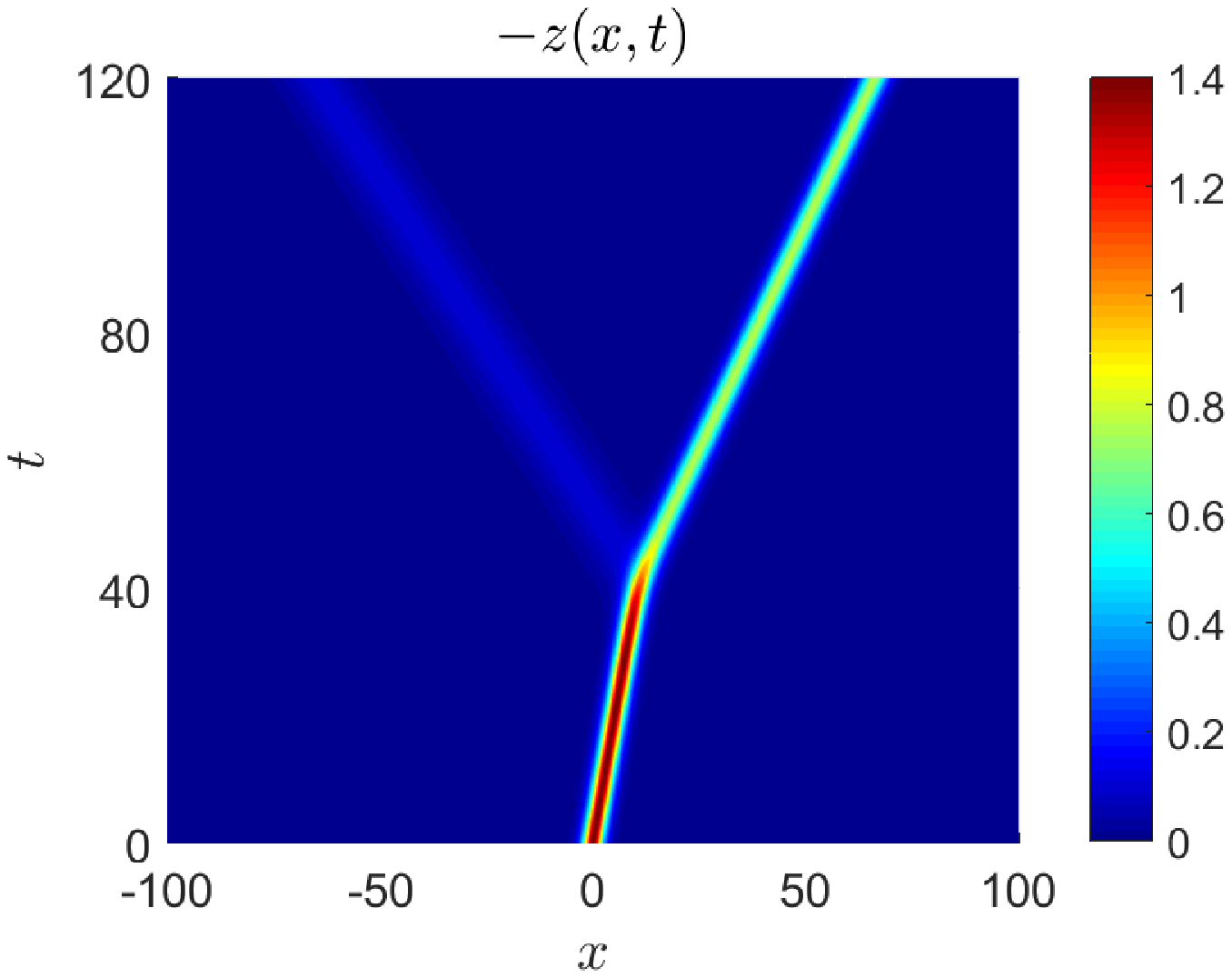}
\end{minipage}%
\hspace{3mm}
\begin{minipage}[t]{0.5\linewidth}
\centering
\includegraphics[height=3cm,width=7.1cm]{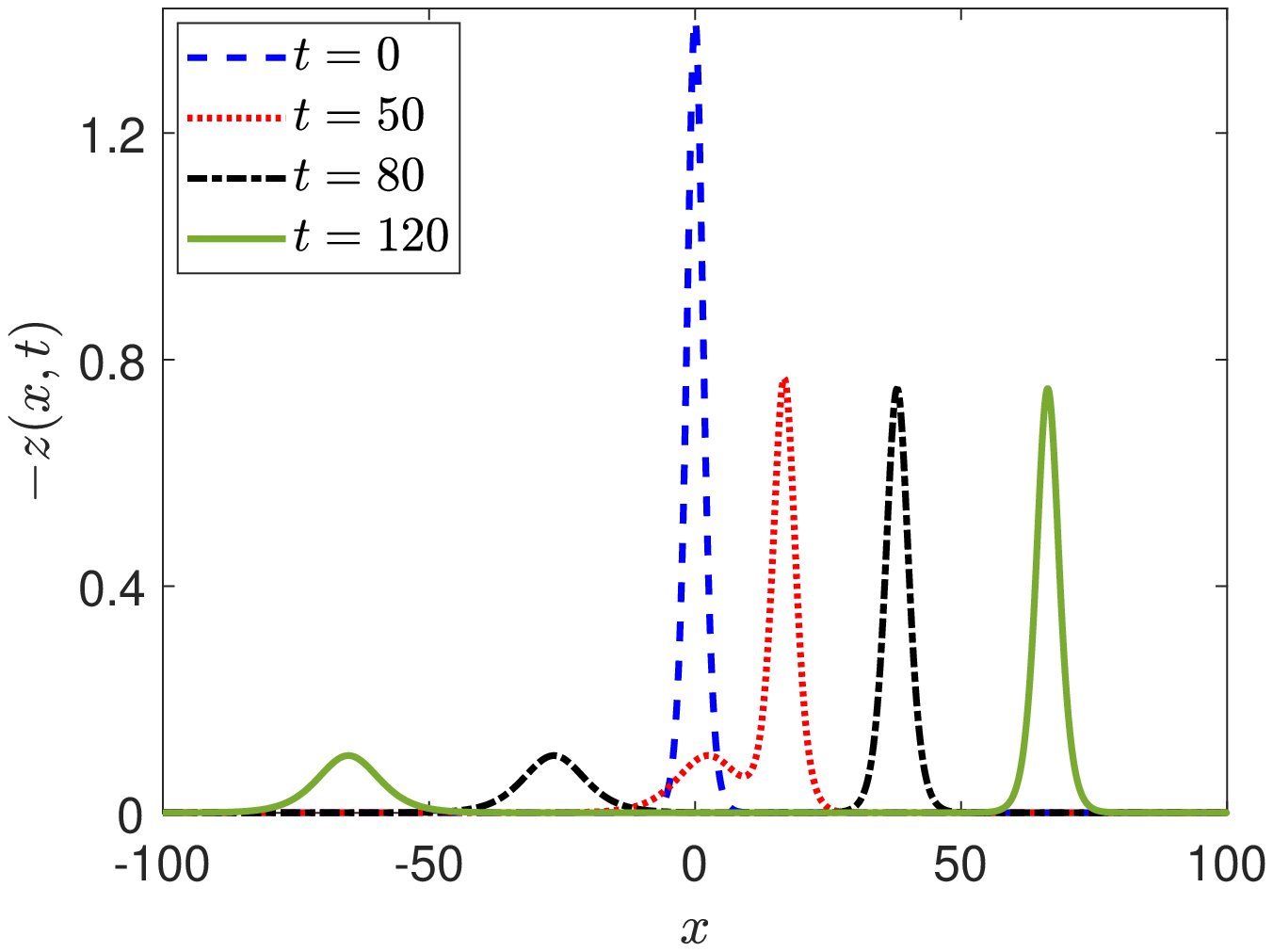}
\end{minipage}
\begin{minipage}[t]{0.5\linewidth}
\centering
\includegraphics[height=3cm,width=7.1cm]{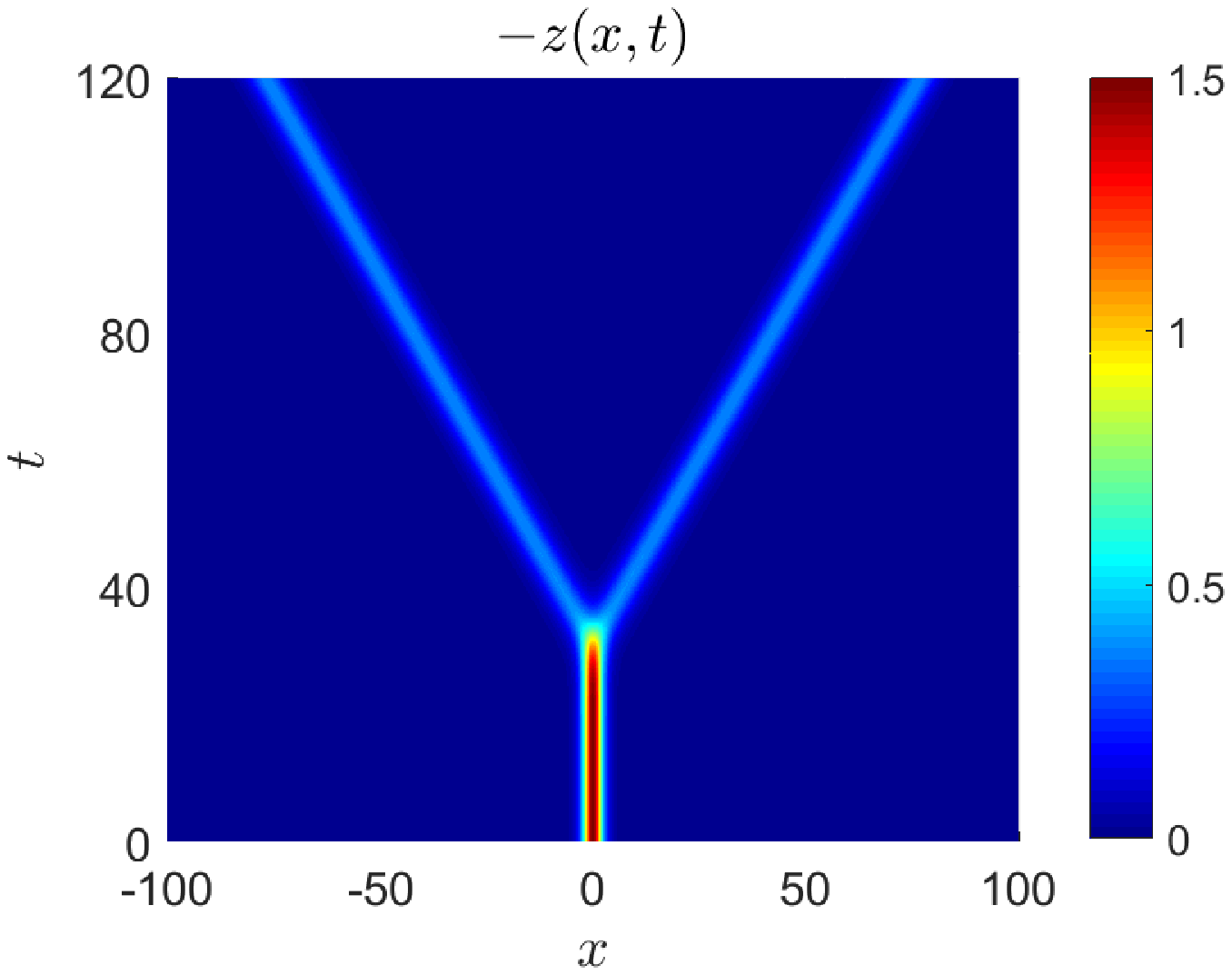}
\end{minipage}%
\hspace{3mm}
\begin{minipage}[t]{0.5\linewidth}
\centering
\includegraphics[height=3cm,width=7.1cm]{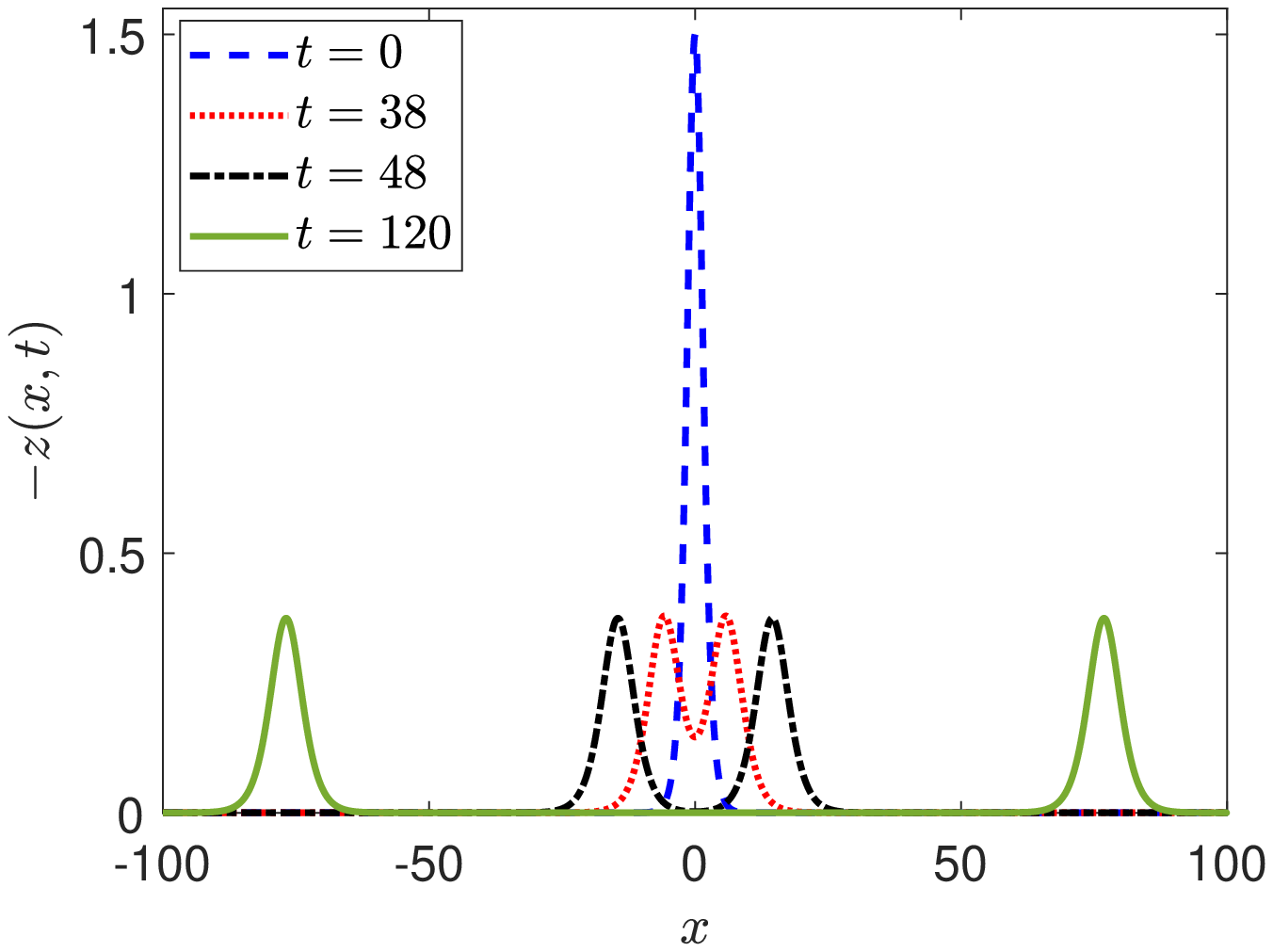}
\end{minipage}
\vspace{-3mm}
\caption{Evolution of the solitary wave in Example 2 for different $A$ :  $A=1.1, 1.3, 1.4, 1.5$ (from top to bottom).}\label{birth}
\end{figure}

Computations are done for $h=1/8$ and $\tau=0.001$ on the interval $\Omega=[-400, 400]$. We find that for small $A$, the numerical solution agrees with the exact solution well. However, for larger $A$, the soliton fails to preserve its shape and splits into two pulses as time evolves. Fig. \ref{birth} shows the evolution of the soliton for different amplitudes $A=1.1, 1.3, 1.4, 1.5$. It can be observed that for smaller $A$, the soliton preserves its shape and velocity well, which agrees with the exact solitary wave solution \eqref{sol-ex}. However, for the initial soliton with larger amplitude, the initial pulse splits into two solitons moving in the opposite directions. Furthermore, the larger $A$ is, the smaller the difference between the amplitudes of the two solitons is. Particularly, for $A=1.5$ which corresponds to null initial velocity $u_1(x)=0$, the soliton finally splits into two solitons with equal amplitudes and equal velocities traveling in opposite directions.

\begin{figure}[h!]
\begin{minipage}[t]{0.5\linewidth}
\centering
\includegraphics[height=3.8cm,width=7.1cm]{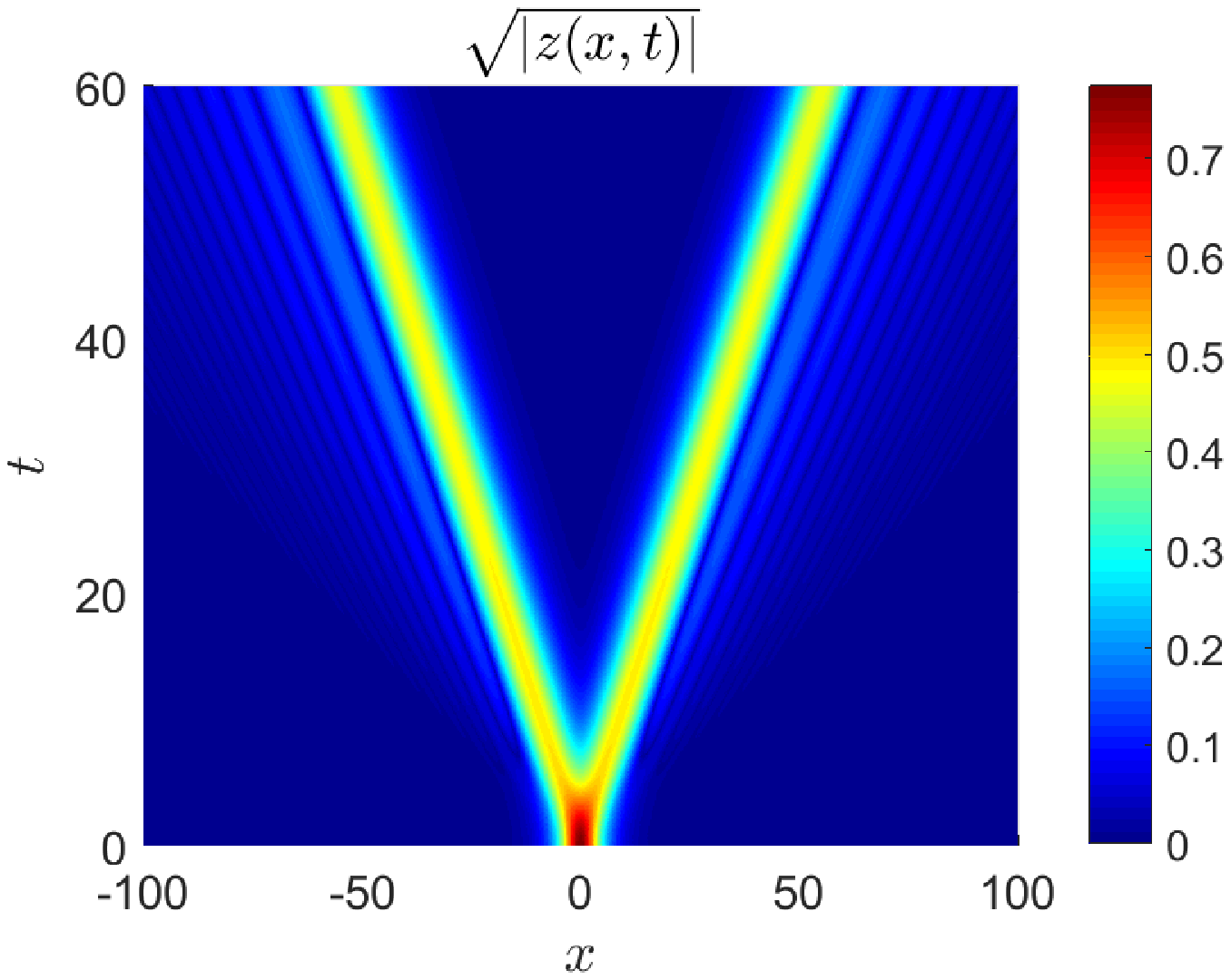}
\end{minipage}%
\hspace{3mm}
\begin{minipage}[t]{0.5\linewidth}
\centering
\includegraphics[height=3.8cm,width=7.1cm]{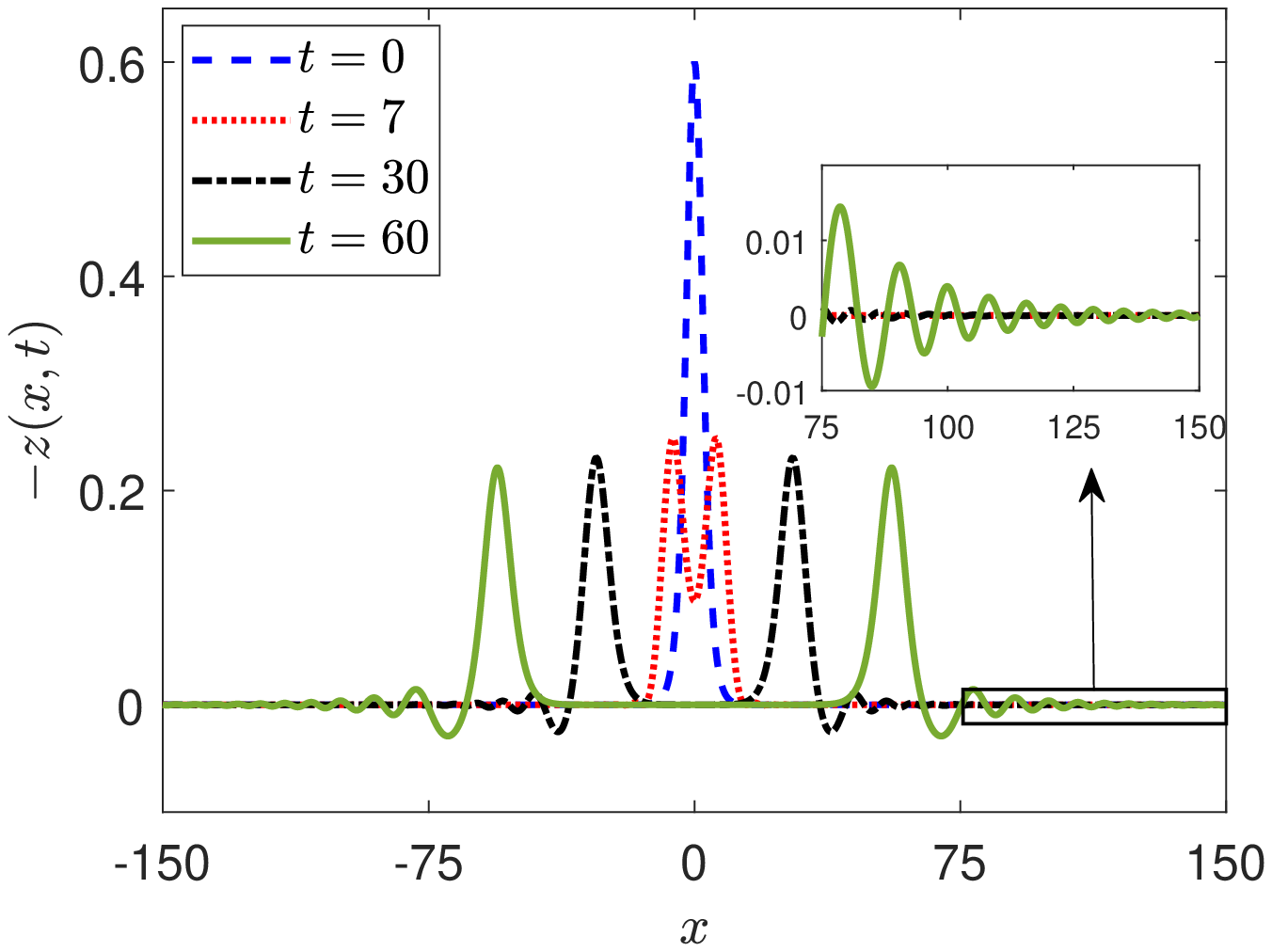}
\end{minipage}
\begin{minipage}[t]{0.5\linewidth}
\centering
\includegraphics[height=3.8cm,width=7.1cm]{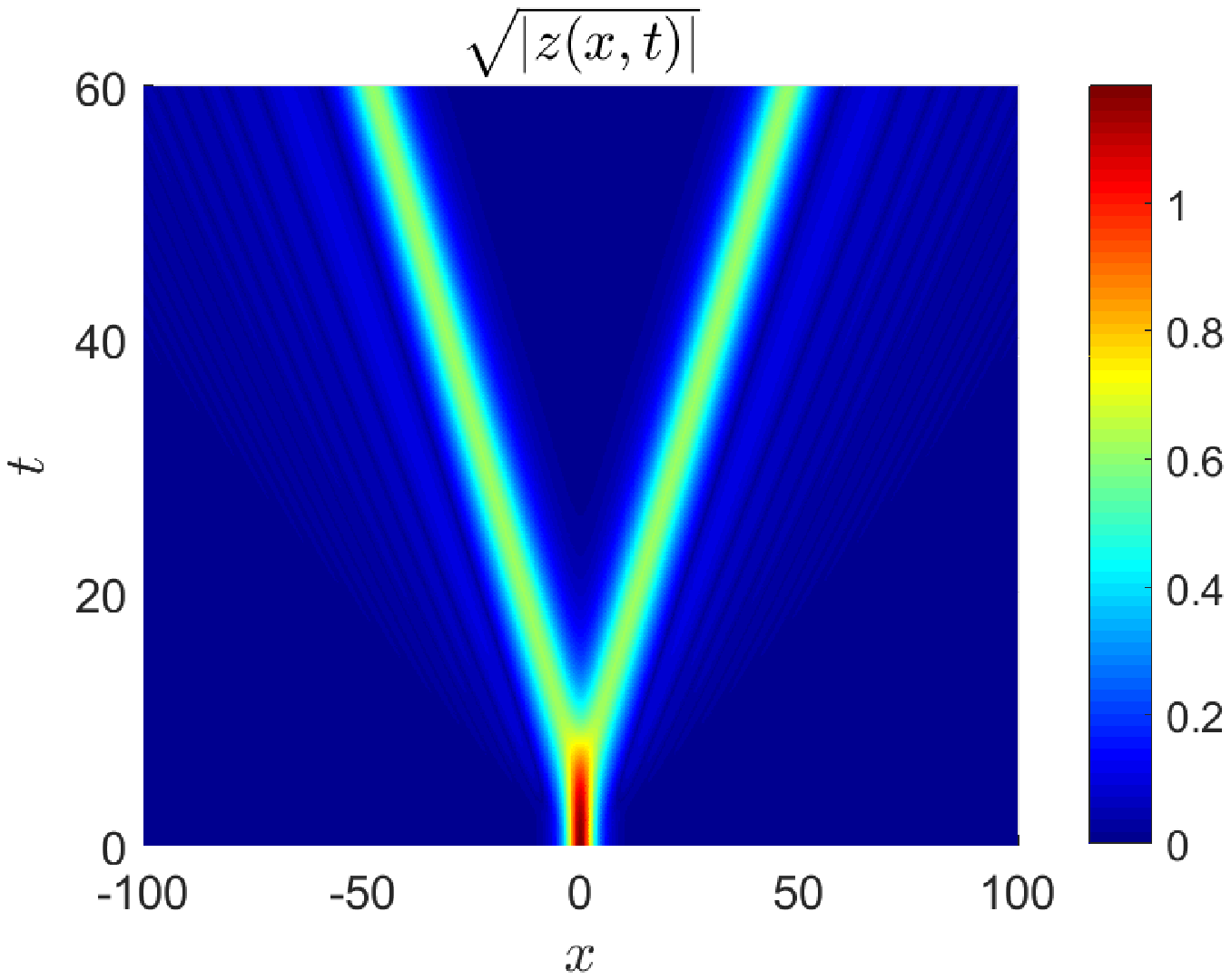}
\end{minipage}%
\hspace{3mm}
\begin{minipage}[t]{0.5\linewidth}
\centering
\includegraphics[height=3.8cm,width=7.1cm]{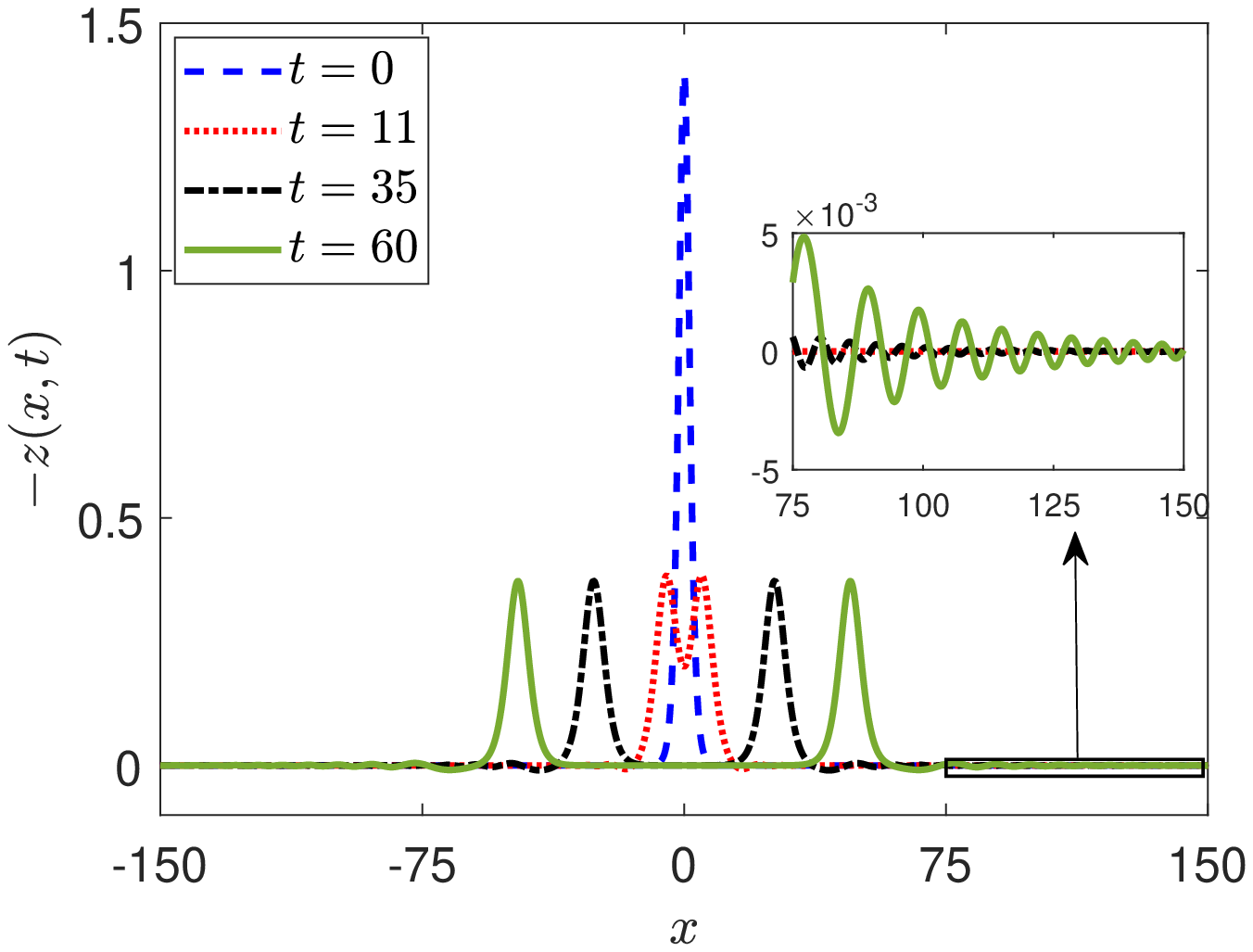}
\end{minipage}
\vspace{-3mm}
\caption{Evolution of the solitary wave with zero initial velocity for different $A$ :  $A=0.6, 1.4$ (from top to bottom).}\label{birth1}
\end{figure}

For comparison, we also investigate the evolution of arbitrary pulse with zero initial velocity:
\be\label{ini_set21}
u_0(x)=-A\,\sech^2(\sqrt{A/6}\,x),\quad u_1(x)=0,
\ee
which has been studied in the literature \cite{jiang2016}. Fig. \ref{birth1} displays the dynamics of the soliton with null initial velocity and different amplitudes $A=0.6, 1.4$.  Different from the case with non-zero initial velocity, the soliton always splits into two pulses with equal amplitudes and equal velocities propagating in opposite directions. Besides the two main solitons, some dispersive oscillations are emitted after splitting as long as $A<1.5$.  Specifically, the larger the amplitude is, the stronger the additional emitting wave is.
\subsection{Interaction of two solitons}
In this section, we apply the DEI-FP method to investigate the interaction of two solitary waves traveling in the opposite directions. The initial data is chosen as
\be\label{ini_set3}
\begin{split}
&u_0(x)=-\sum_{k=1}^2  A_k\,\sech^2\left(\sqrt{\fl{A_k}{6}}(x-x_k)\right),\quad v_k=\pm \sqrt{1-\fl{2}{3}A_k},\\
&u_1(x)=-\sum_{k=1}^2  A_k v_k\sqrt{\fl{2A_k}{3}}
\,\sech^2\left(\sqrt{\fl{A_k}{6}}(x-x_k)\right)\tanh\left(\sqrt{\fl{A_k}{6}}(x-x_k)\right).
\end{split}
\ee
It represents two solitary waves located initially at the positions $x=x_1$ and $x=x_2$, respectively, moving to the right or left depending on the sign of the velocity $v_k$. Computations are done for $h=1/8$ and $\tau=0.001$ on the interval $\Omega=[-400, 400]$. We consider the following cases:

\noindent(1) \emph{Elastic collision:}
\begin{itemize}
\item[](i).    $x_1=-x_2=50$, $A_1=0.2$, $A_2=0.3$, $v_1>0$, $v_2<0$;
\item[] (ii).   $x_1=-x_2=10$,  $A_1=0.2$, $A_2=0.5$, $v_1>0$, $v_2<0$;
\end{itemize}

\noindent(2) \emph{Blow-up phenomenon:}
\begin{itemize}
\item[] (iii). $x_1=-x_2=50$,  $A_1=0.37$, $A_2=0.37$, $v_1>0$, $v_2<0$;
\item[] (iv). $x_1=-x_2=50$,  $A_1=0.38$, $A_2=0.38$, $v_1>0$, $v_2<0$;
\item[] (v).  $x_1=-x_2=50$,  $A_1=0.3$, $A_2=0.45$, $v_1>0$, $v_2<0$;
\item[] (vi). $x_1=-x_2=50$,  $A_1=0.3$, $A_2=0.46$, $v_1>0$, $v_2<0$;
\end{itemize}

\noindent(3) \emph{Interaction with static solitons:}
\begin{itemize}
\item[] (vii). $x_1=-x_2=50$,  $A_1=0.37$, $A_2=1.5$, $v_1>0$;
\item[] (viii).  $x_1=-x_2=50$,  $A_1=0.38$, $A_2=1.5$, $v_1>0$;
\item[] (ix). $x_1=-x_2=30$,  $A_1=A_2=1.5$, $v_1=v_2=0$;
\item[] (x).  $x_1=-x_2=20$,  $A_1=A_2=1.5$, $v_1=v_2=0$;
\end{itemize}

\noindent(4) \emph{Overtaking interaction:}
\begin{itemize}
\item[] (xi).  $x_1=-80$, $x_2=-50$,  $A_1=0.2$, $A_2=1$, $v_1>0$, $v_2>0$.
\end{itemize}

\begin{figure}[h!]
\begin{minipage}[t]{0.5\linewidth}
\centering
\includegraphics[height=3.5cm,width=7.1cm]{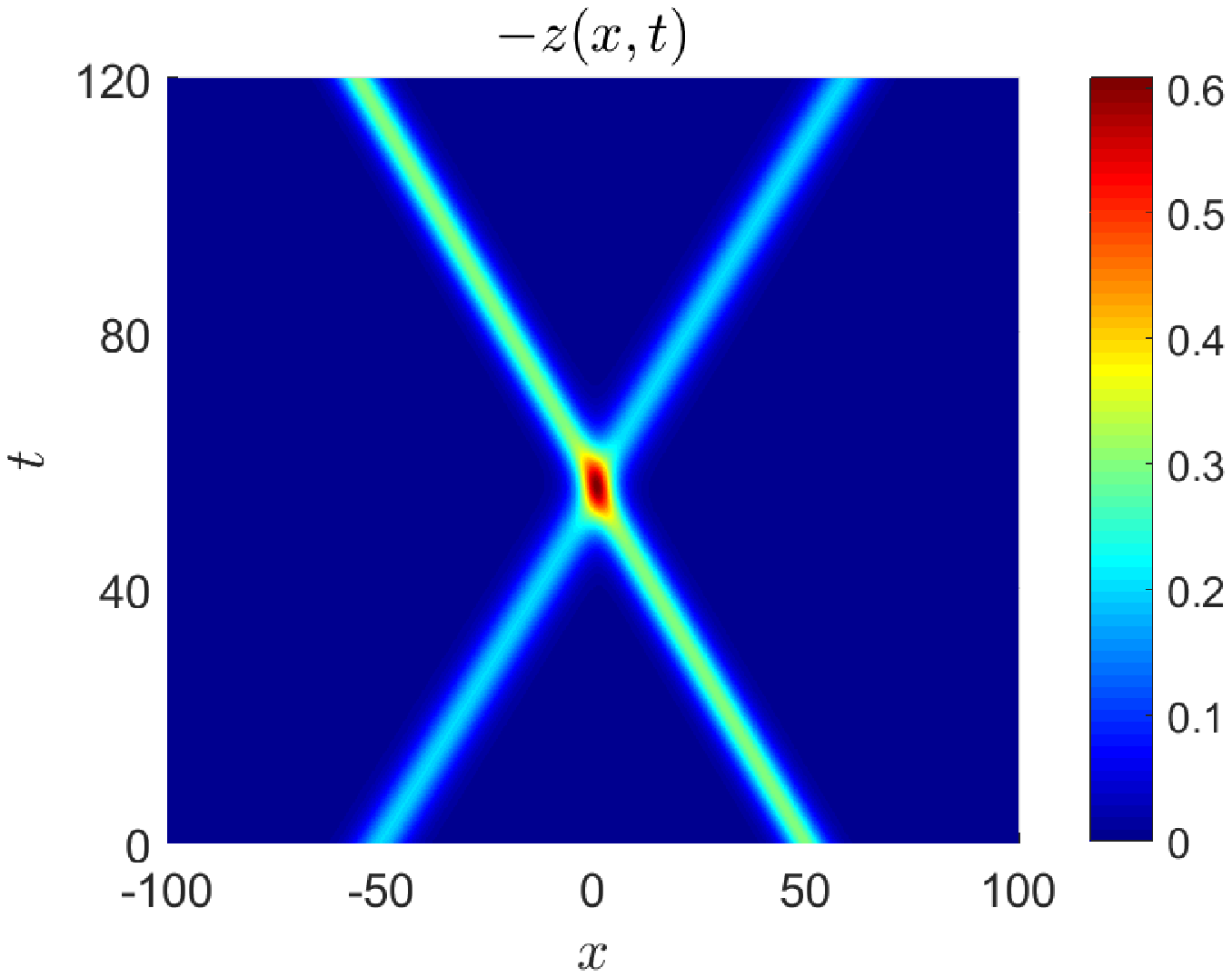}
\end{minipage}%
\hspace{3mm}
\begin{minipage}[t]{0.5\linewidth}
\centering
\includegraphics[height=3.5cm,width=7.1cm]{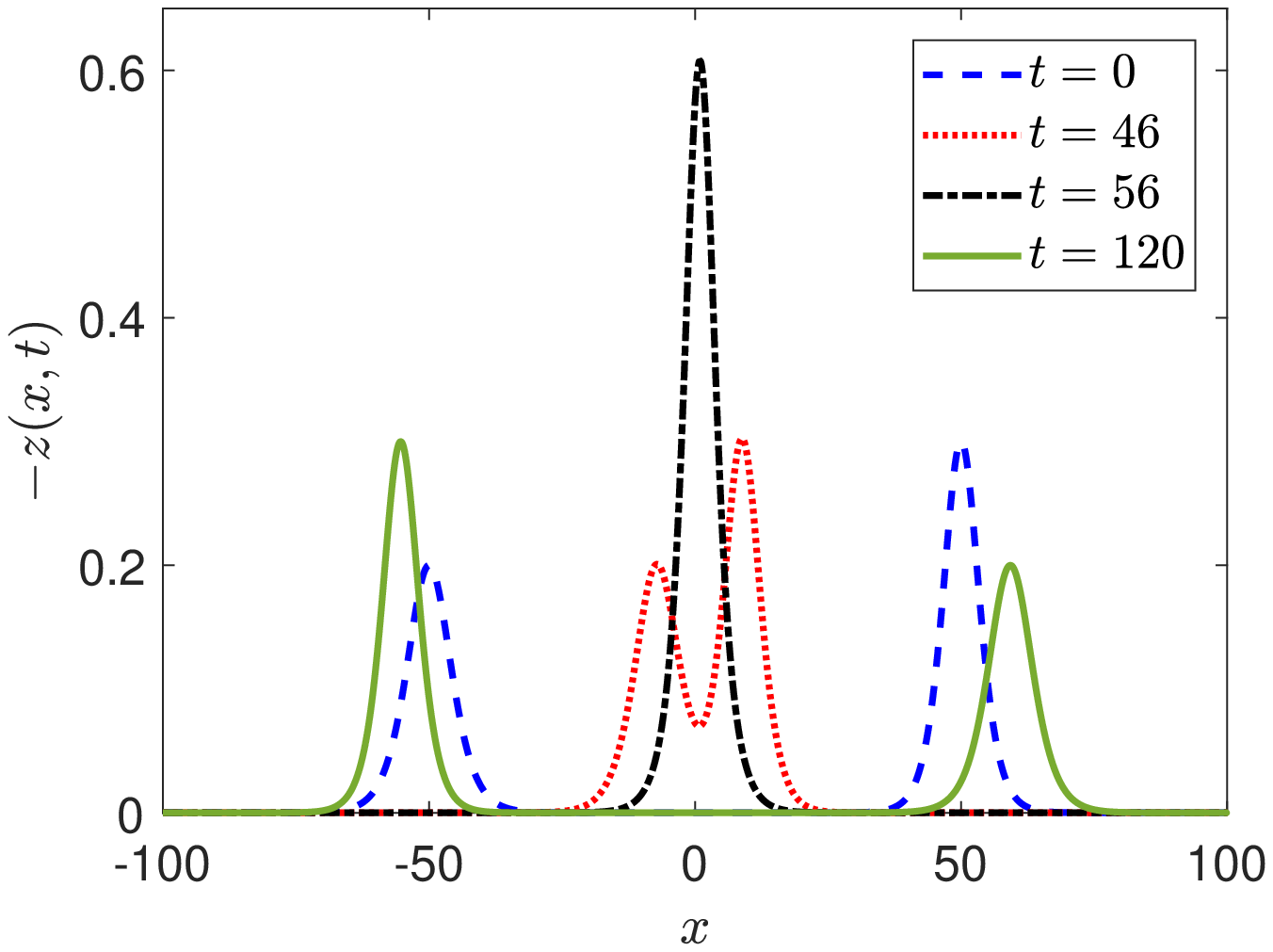}
\end{minipage}
\begin{minipage}[t]{0.5\linewidth}
\centering
\includegraphics[height=3.5cm,width=7.1cm]{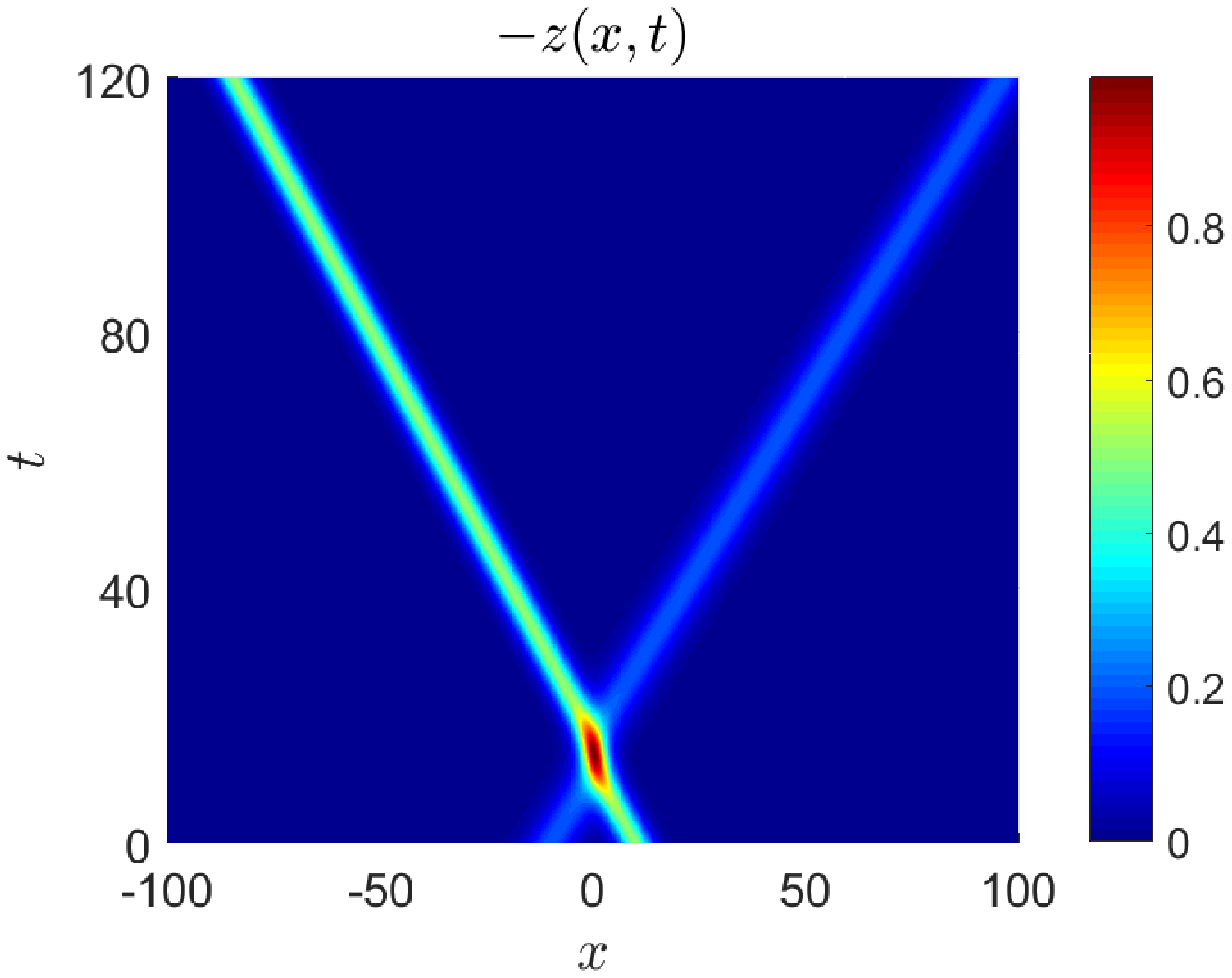}
\end{minipage}%
\hspace{3mm}
\begin{minipage}[t]{0.5\linewidth}
\centering
\includegraphics[height=3.5cm,width=7.1cm]{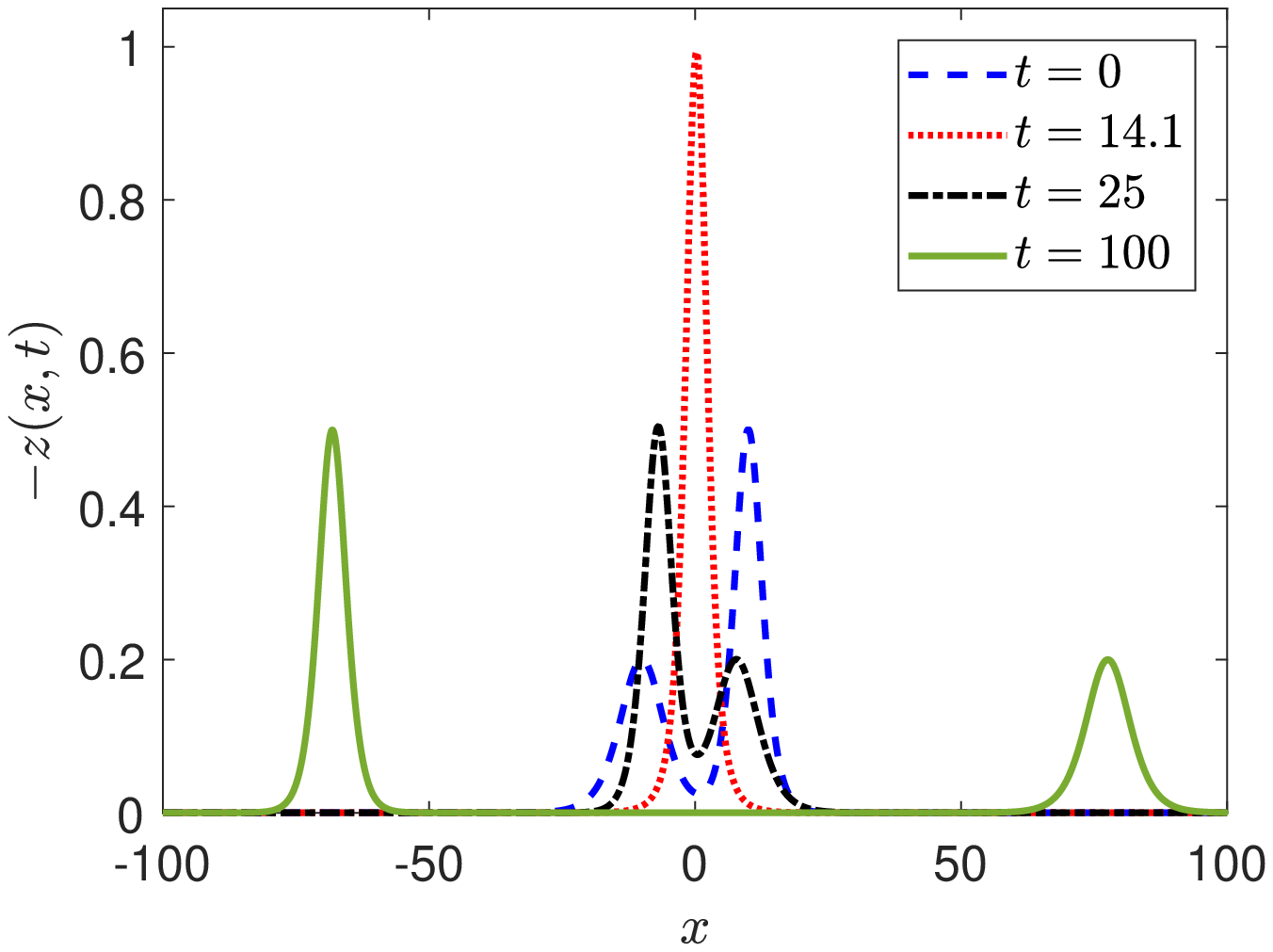}
\end{minipage}
\vspace{-3mm}
\caption{Elastic collision of two solitons for Cases (i)-(ii) (from top to bottom).}\label{interac}
\end{figure}

\begin{figure}[h!]
\begin{minipage}[t]{0.5\linewidth}
\centering
\includegraphics[height=3.5cm,width=7.1cm]{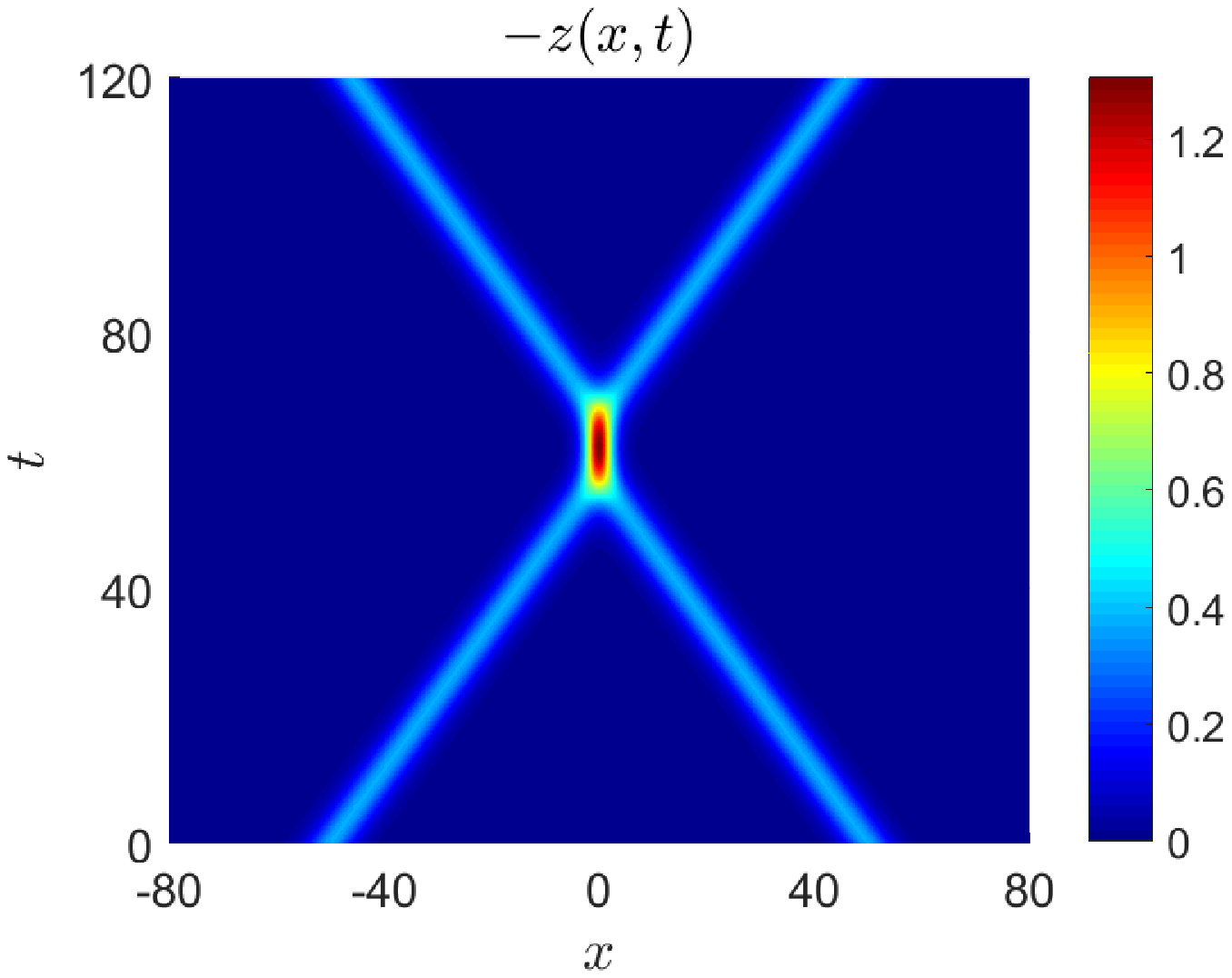}
\end{minipage}%
\hspace{3mm}
\begin{minipage}[t]{0.5\linewidth}
\centering
\includegraphics[height=3.5cm,width=7.1cm]{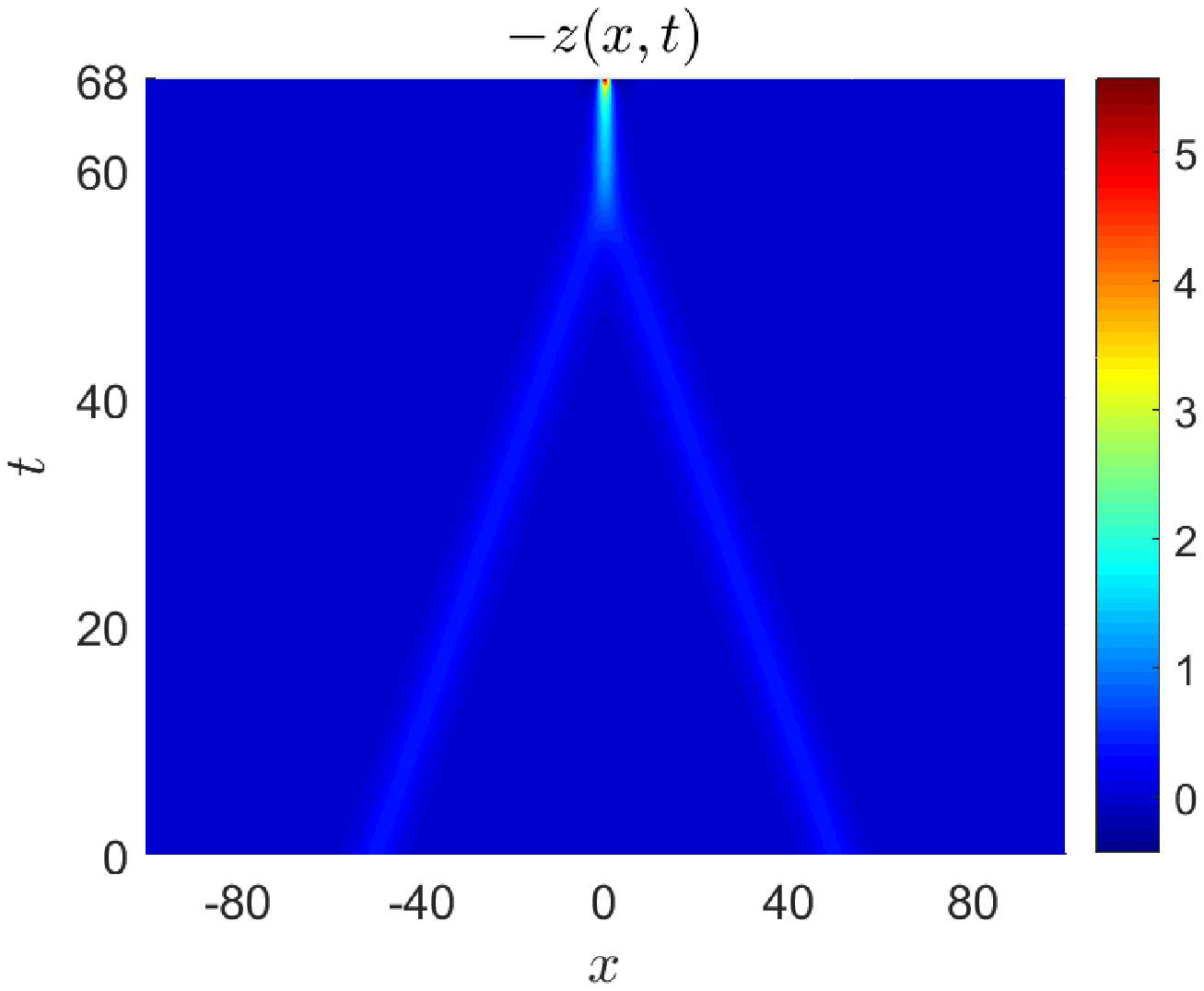}
\end{minipage}
\begin{minipage}[t]{0.5\linewidth}
\centering
\includegraphics[height=3.5cm,width=7.1cm]{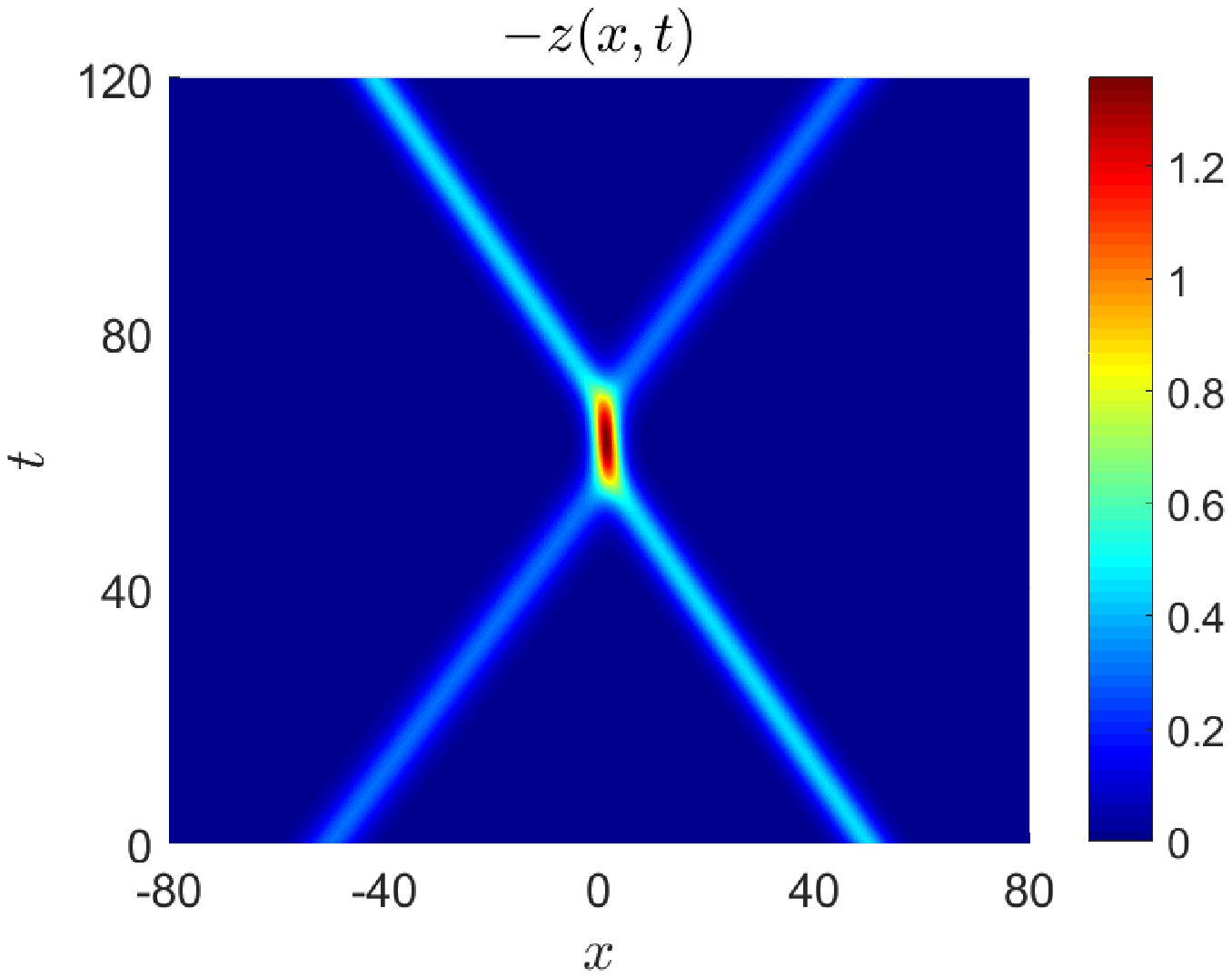}
\end{minipage}%
\hspace{3mm}
\begin{minipage}[t]{0.5\linewidth}
\centering
\includegraphics[height=3.5cm,width=7.1cm]{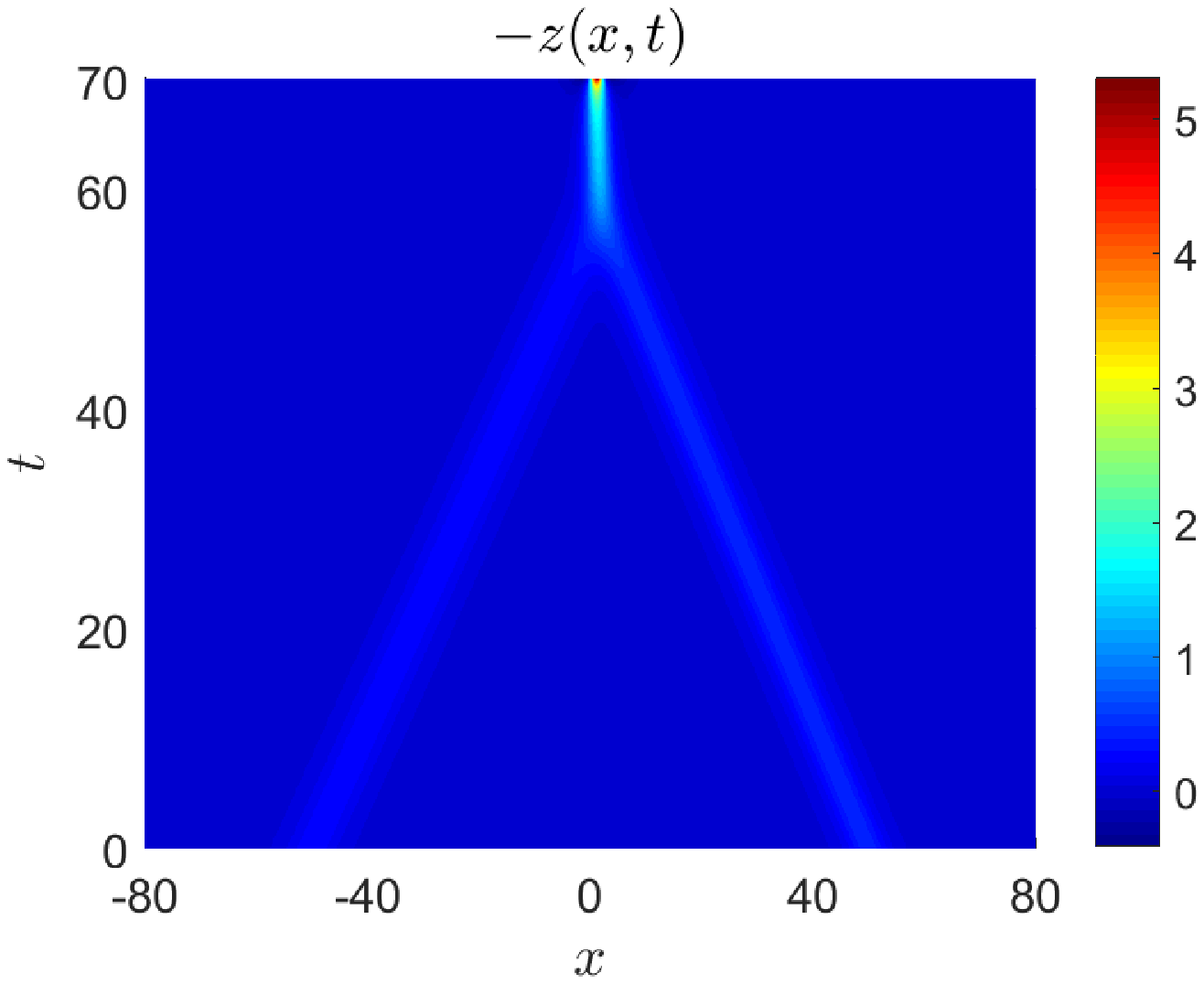}
\end{minipage}
\vspace{-3mm}
\caption{Collision of two solitons for Cases (iii)-(vi) (from top to bottom, from left to right).}\label{blow-up}
\end{figure}

Fig. \ref{interac} shows the evolution of $-z(x,t)$ at different time for \emph{elastic collision} (Cases (i)-(ii)). We see that the two solitons which are initially located at the positions $x_1=-50$ and $x_2=50$ moving towards each other with velocities $v_1$ and $v_2$, respectively. As time progresses they collide, stick together and split after collision without changing their shape and velocities. It can be clearly seen that the collision  generates no radiation. Similar phenomena occurs for Case (ii) where the two solitons contact each other initially.

Fig. \ref{blow-up} investigates the \emph{blow-up phenomenon} for the head-on collision. For $A_1=A_2=0.38$, the solution blows up quickly after $t=68$. Similar blow-up occurs for $A_1=0.3$, $A_2=0.46$ after $t=70$. We see that for $A_1=A_2$, there exists $A_c\in (0.37, 0.38)$ such that the solution blows up in finite time when $A_1=A_2>A_c$. Similarly, for fixed $A_1=0.3$, there exists $A_c\in (0.45, 0.46)$ such that the solution blows up in finite time when $A_2>A_c$. This blow-up phenomenon was revealed in \cite{Ismail} for two solitons with the same initial amplitude.

\begin{figure}[h!]
\begin{minipage}[t]{0.5\linewidth}
\centering
\includegraphics[height=4cm,width=7.1cm]{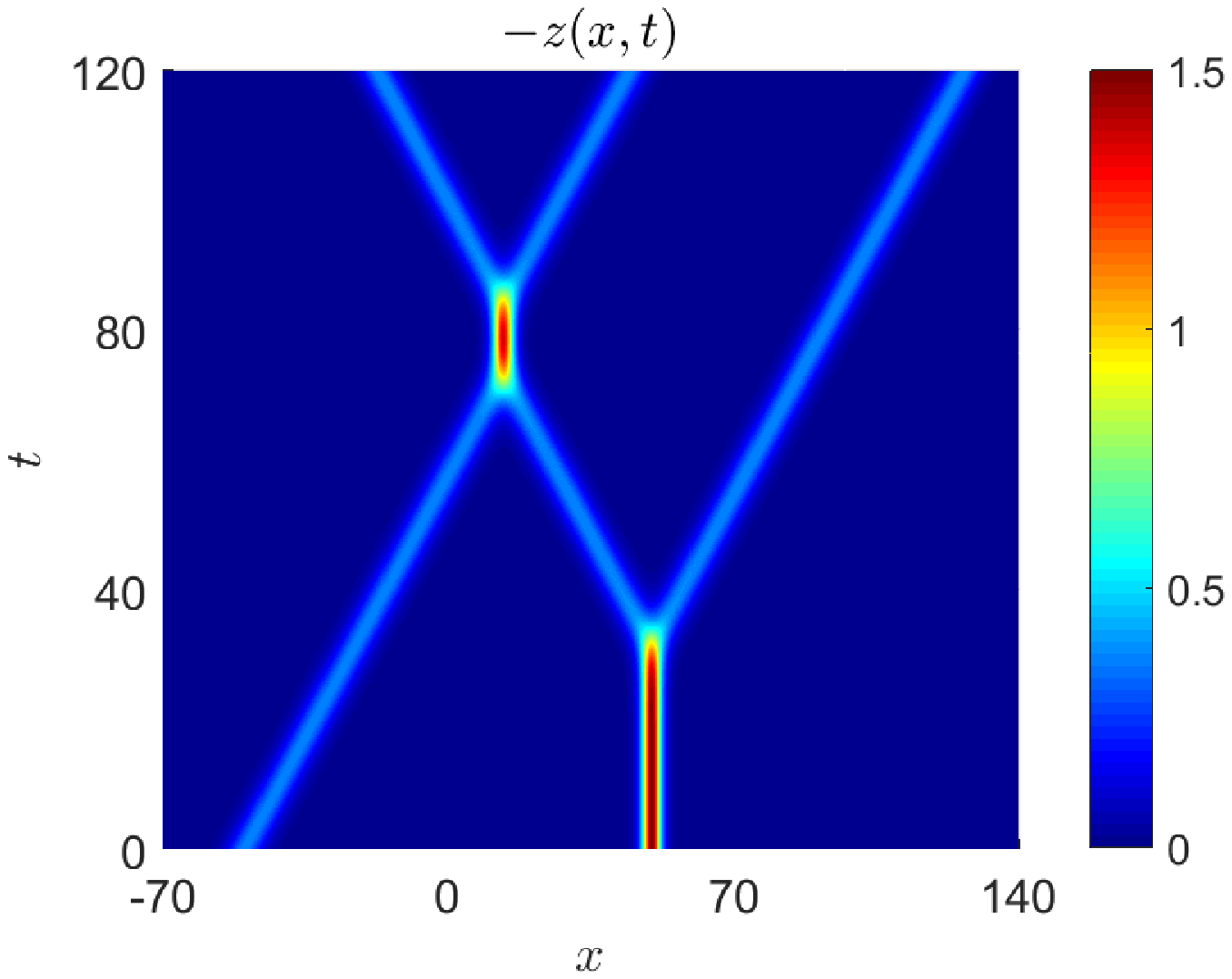}
\end{minipage}%
\hspace{3mm}
\begin{minipage}[t]{0.5\linewidth}
\centering
\includegraphics[height=4cm,width=7.1cm]{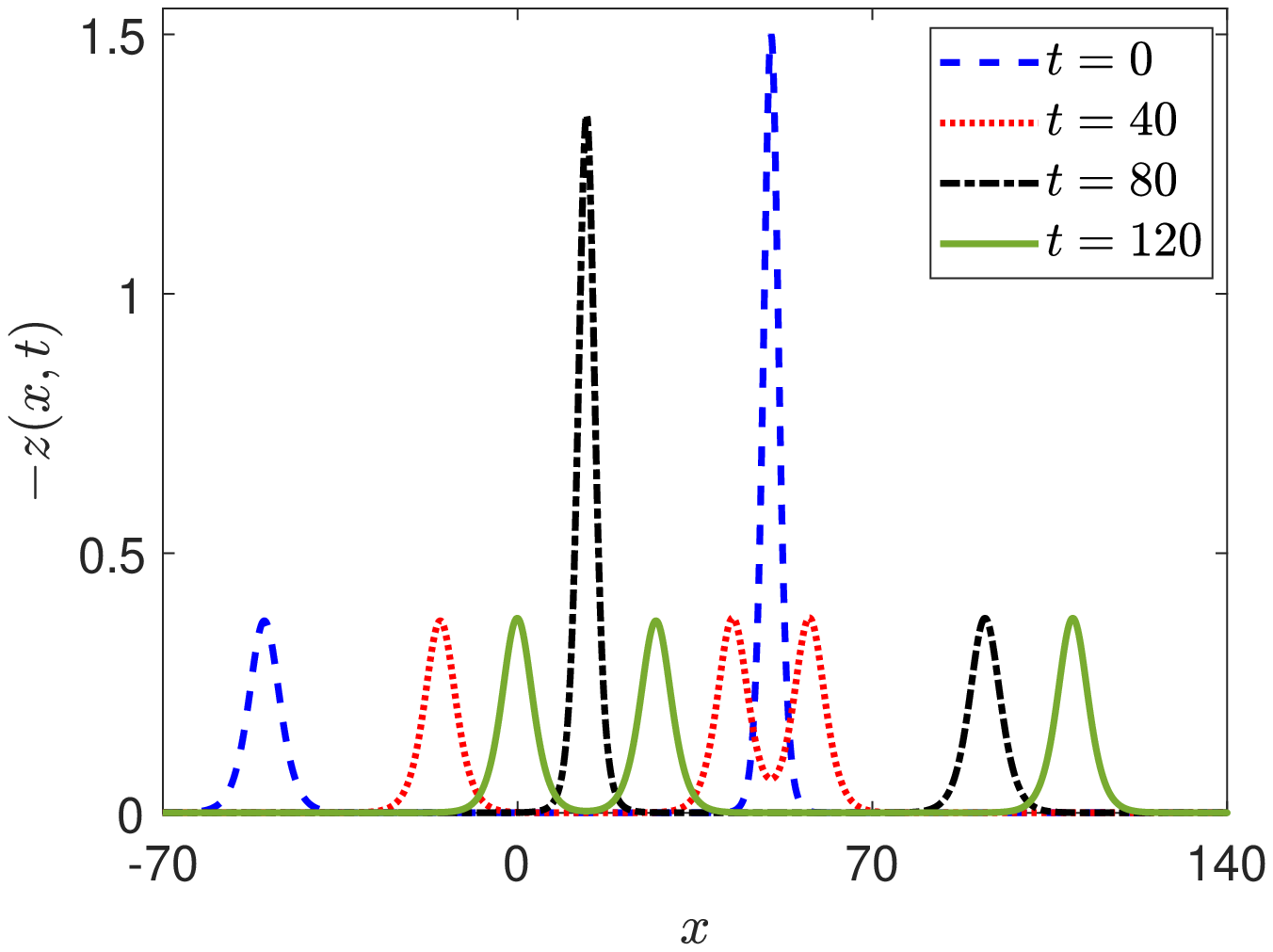}
\end{minipage}
\begin{minipage}[t]{0.5\linewidth}
\centering
\includegraphics[height=4cm,width=7.1cm]{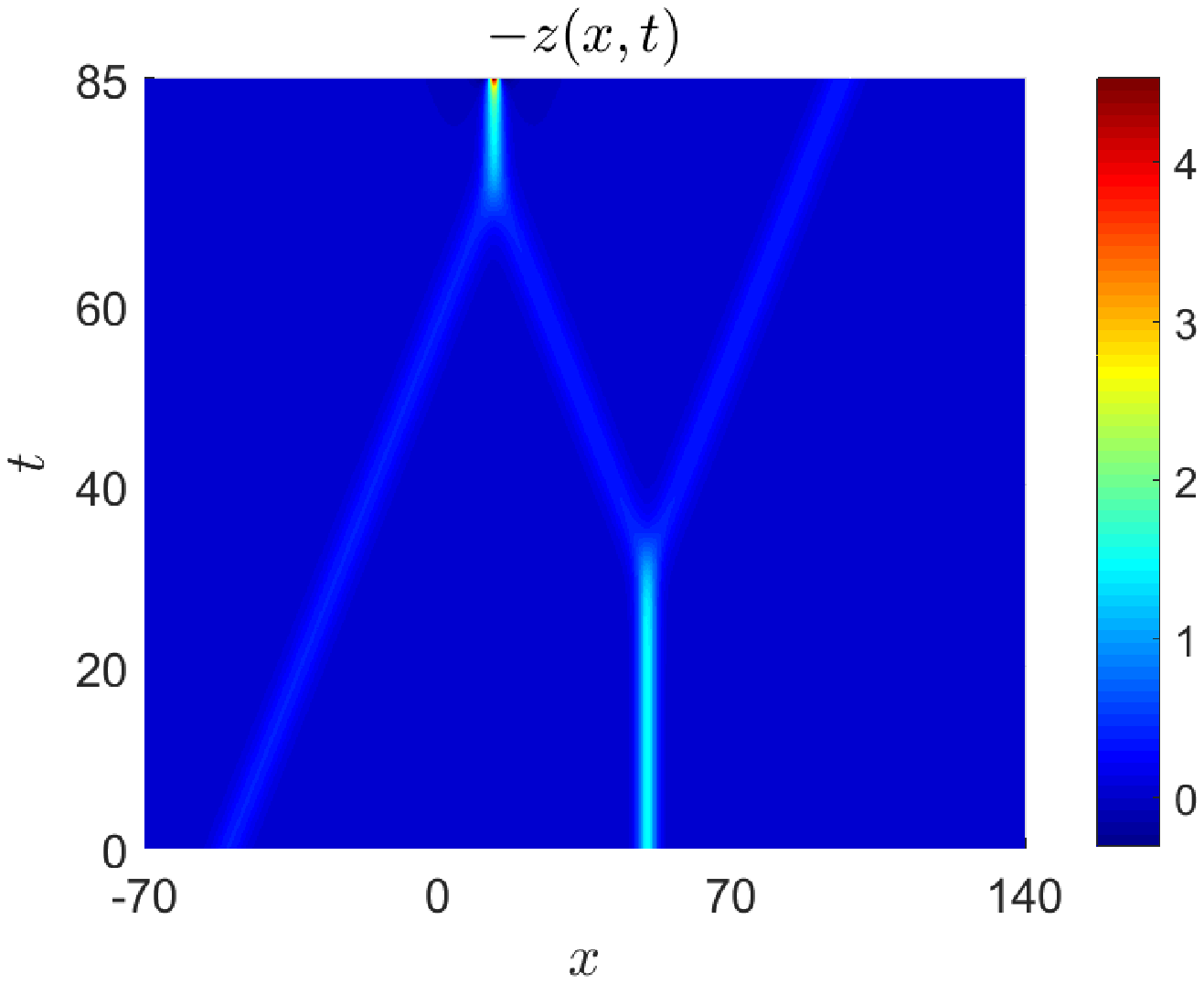}
\end{minipage}%
\hspace{3mm}
\begin{minipage}[t]{0.5\linewidth}
\centering
\includegraphics[height=4cm,width=7.1cm]{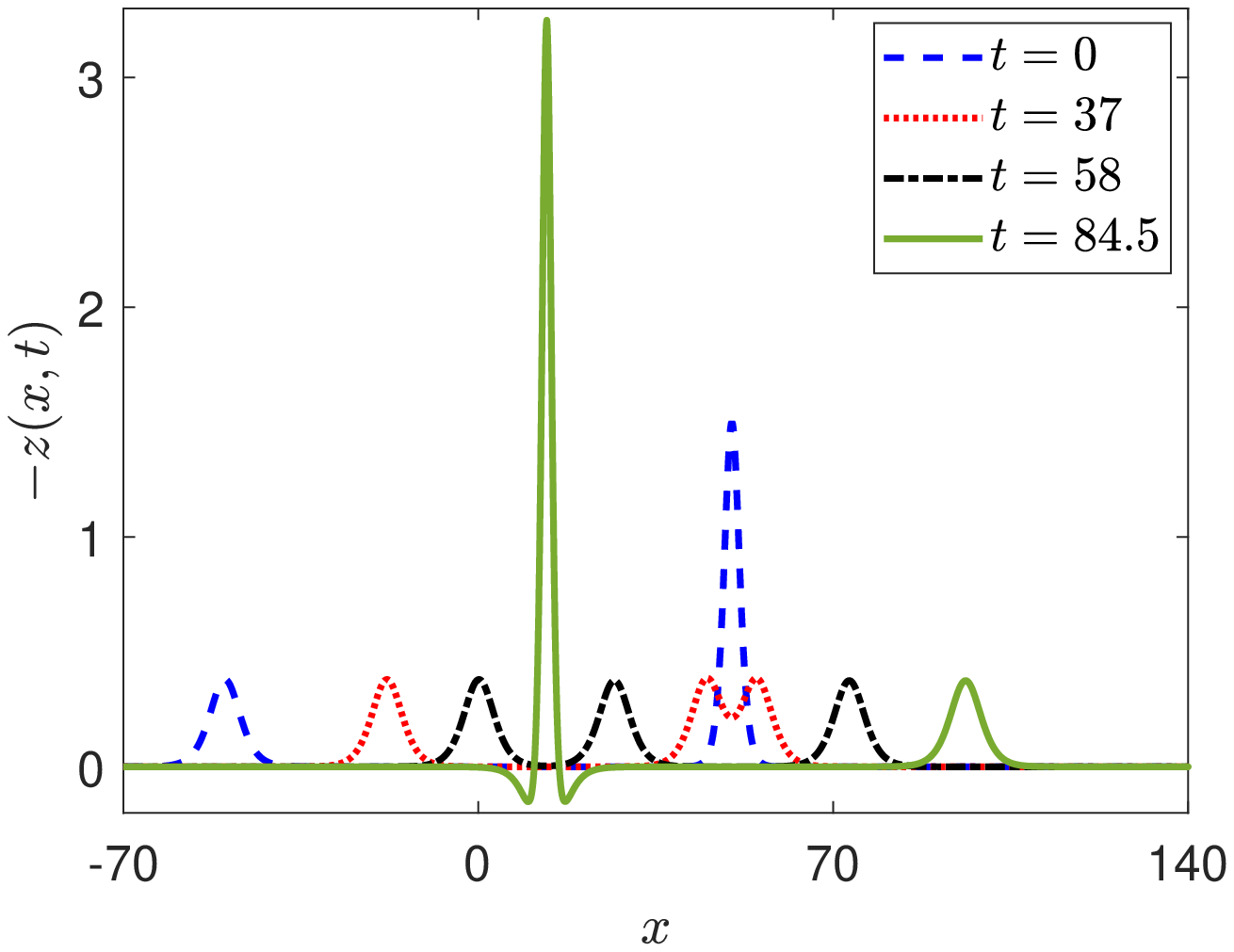}
\end{minipage}
\vspace{-3mm}
\caption{Collision of two solitons for Cases (vii)-(viii) (from top to bottom).}\label{static-in}
\end{figure}

Fig. \ref{static-in} shows the interaction of two solitons, one of which has initial amplitude $A=1.5$ and null initial velocity. As time evolves, the static soliton splits into two pulses propagating in different directions. Then one of the splitting solitons collides with the pulse moving towards it. When the amplitude is not large enough, they collide and split again without changing their shape and velocities (cf. Fig. \ref{static-in} top). While if the amplitude is large enough, blow-up occurs quickly after the collision (cf. Fig. \ref{static-in} bottom).

\begin{figure}[h!]
\begin{minipage}[t]{0.5\linewidth}
\centering
\includegraphics[height=4cm,width=7.1cm]{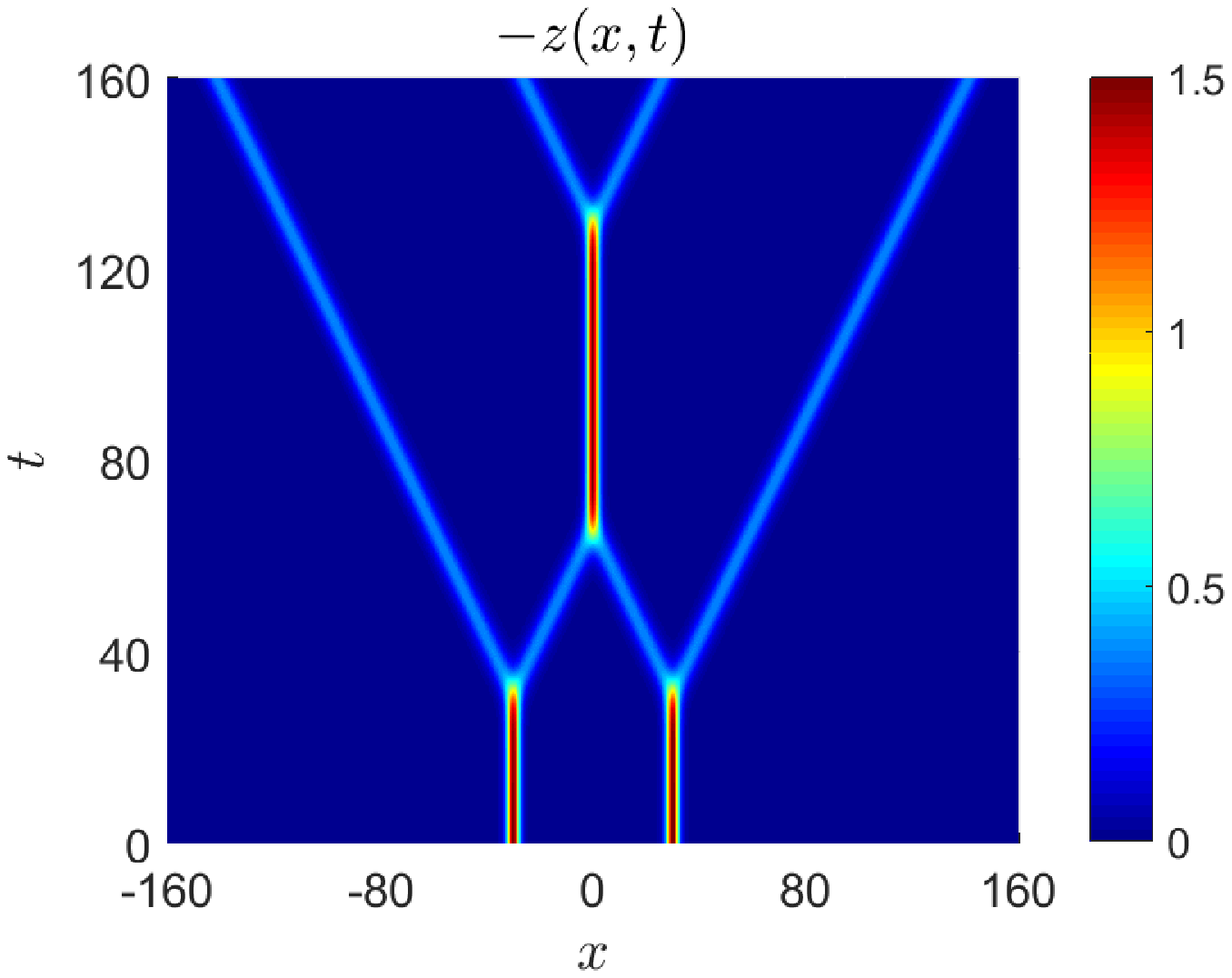}
\end{minipage}%
\hspace{3mm}
\begin{minipage}[t]{0.5\linewidth}
\centering
\includegraphics[height=4cm,width=7.1cm]{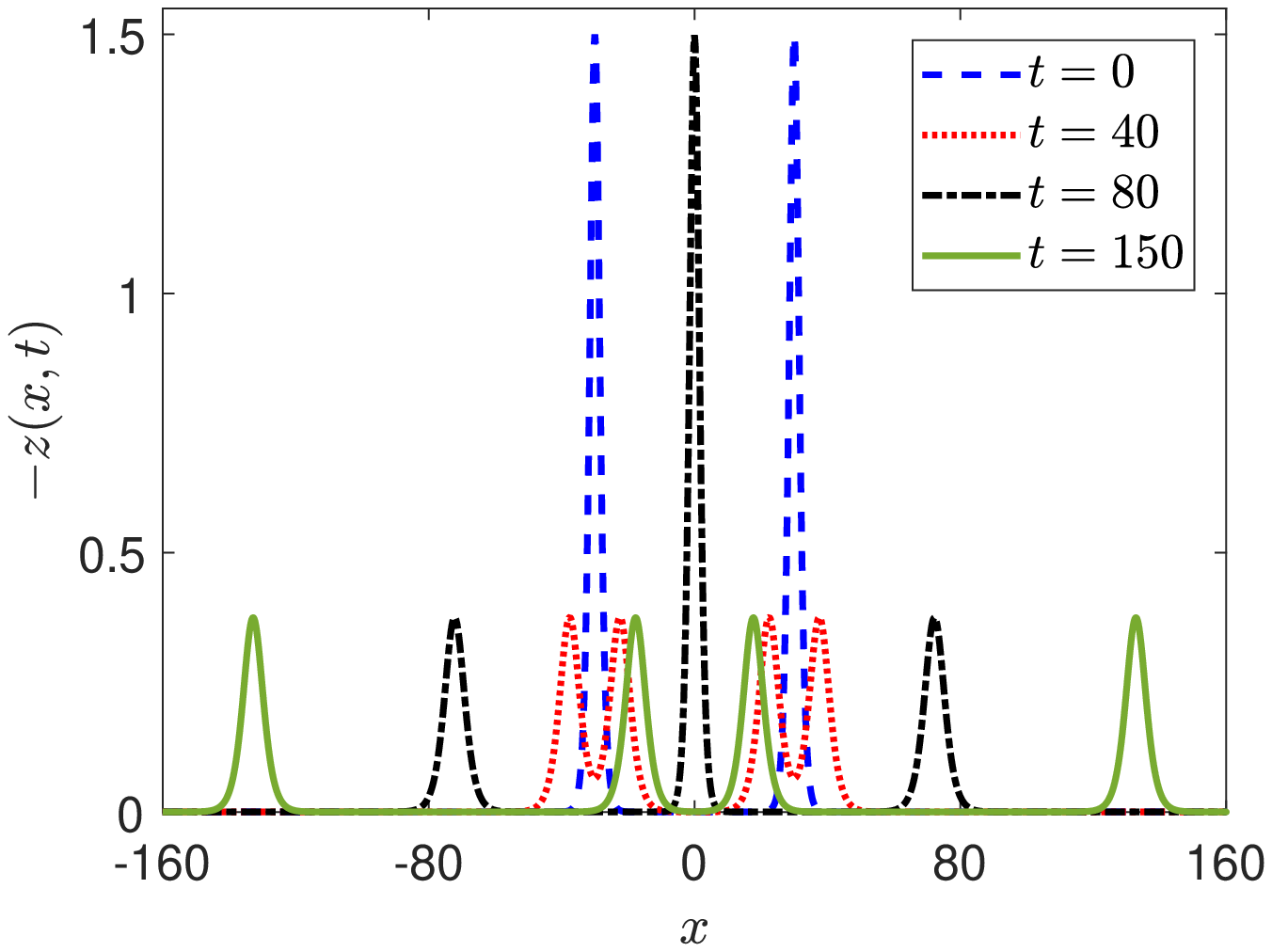}
\end{minipage}
\begin{minipage}[t]{0.5\linewidth}
\centering
\includegraphics[height=4cm,width=7.1cm]{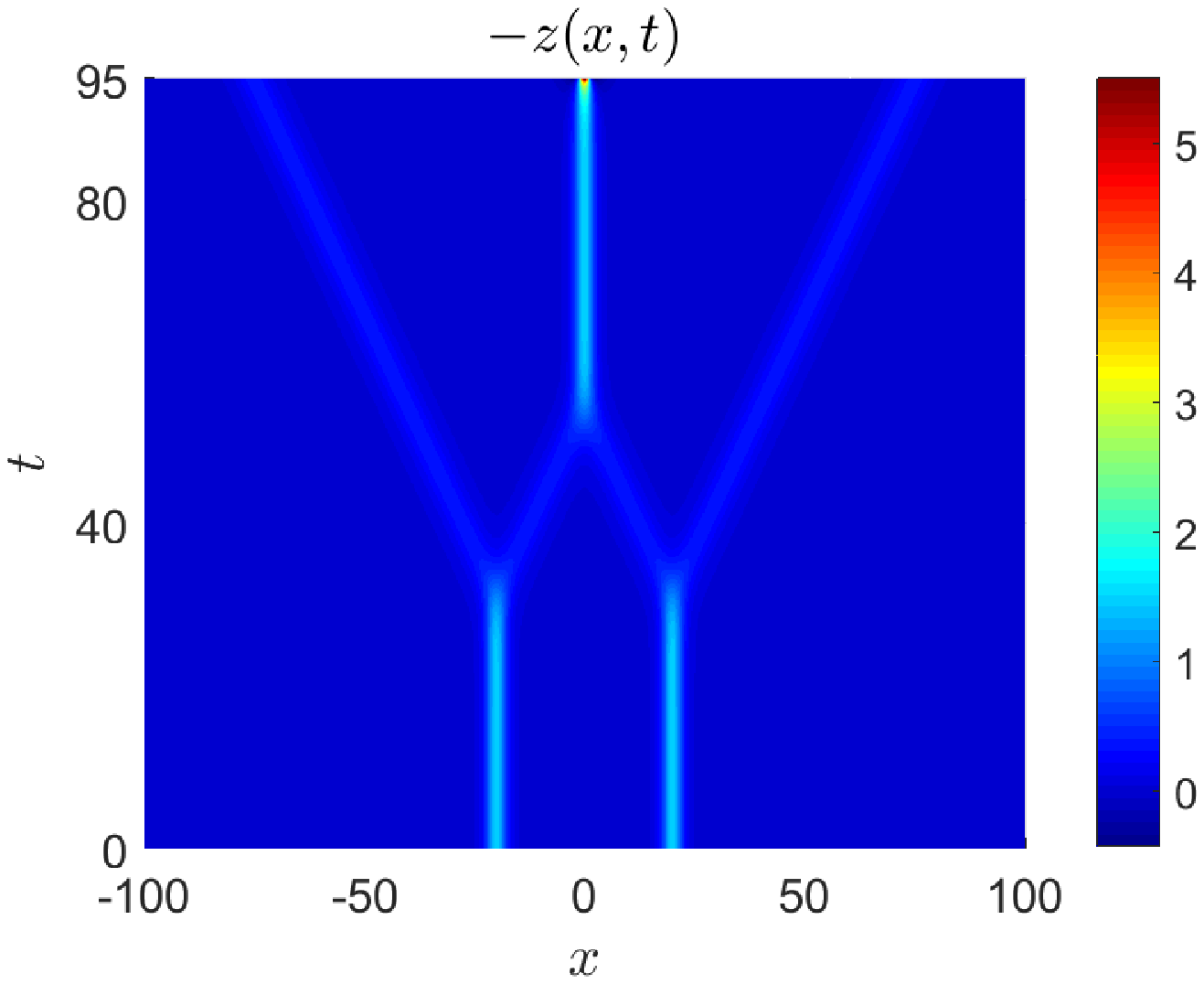}
\end{minipage}%
\hspace{3mm}
\begin{minipage}[t]{0.5\linewidth}
\centering
\includegraphics[height=4cm,width=7.1cm]{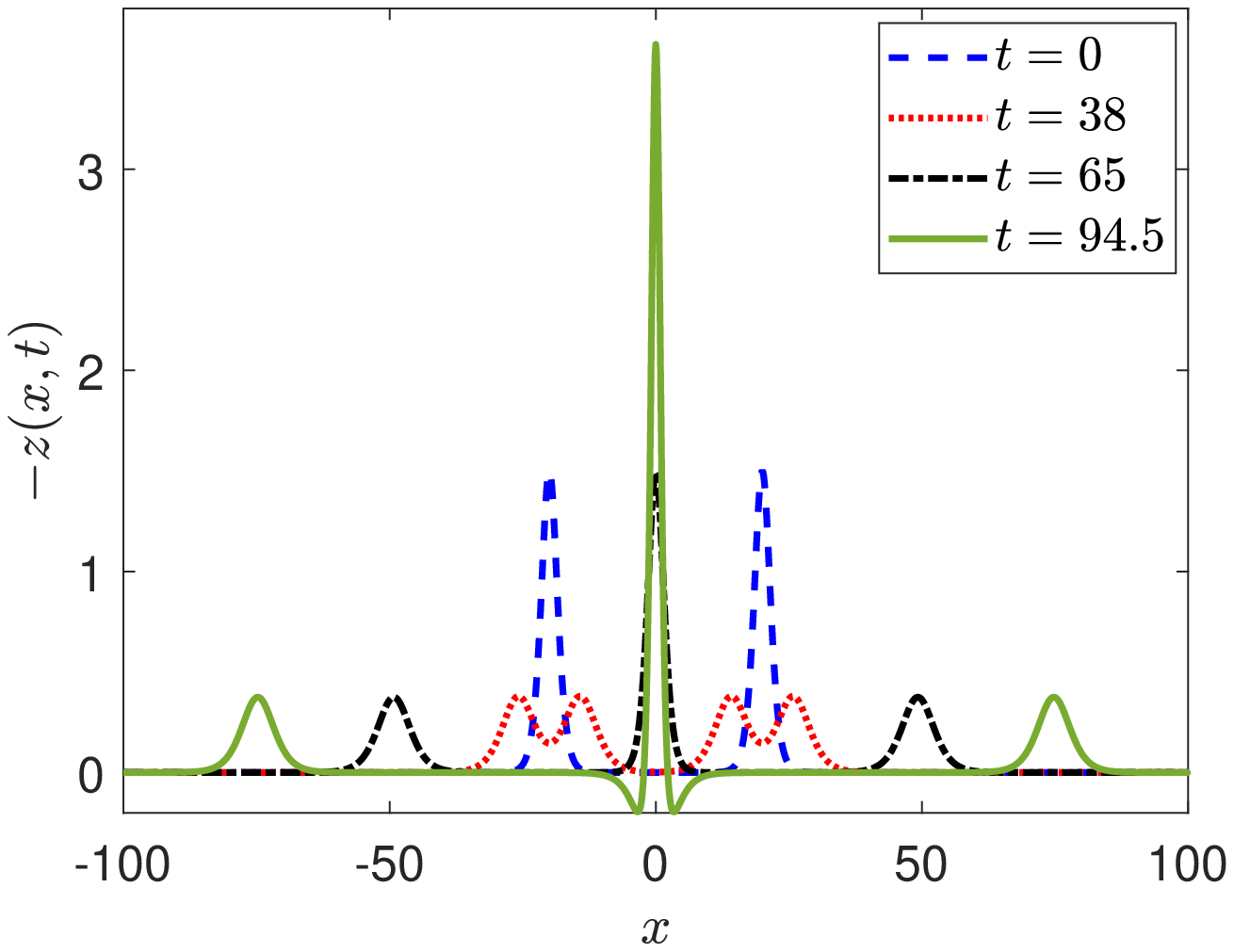}
\end{minipage}
\vspace{-3mm}
\caption{Evolution of two solitons with null initial velocities for Cases (ix)-(x) (from top to bottom).}\label{static2}
\end{figure}

Fig. \ref{static2} displays the evolution of two solitons with initial amplitudes $A_1=A_2=1.5$ and zero initial velocities. When the two solitons are well-separated initially (Case (ix)), each pulse splits into two solitons spreading towards opposite directions. Then the two head-on pulses collide and recover as one solitary wave with amplitude $A=1.5$.  The recovered soliton keeps static for a period of time and finally split into two equal pulses moving in the opposite directions. However, when the two solitons are not initially well-separated (Case (x)), the recovered solitary wave blows up quickly after $t=95$. As far as we know, this phenomena of \emph{recovering as a static soliton} after collision has not been found in the existing literature.

\begin{figure}[h!]
\begin{minipage}[t]{0.5\linewidth}
\centering
\includegraphics[height=4cm,width=7.1cm]{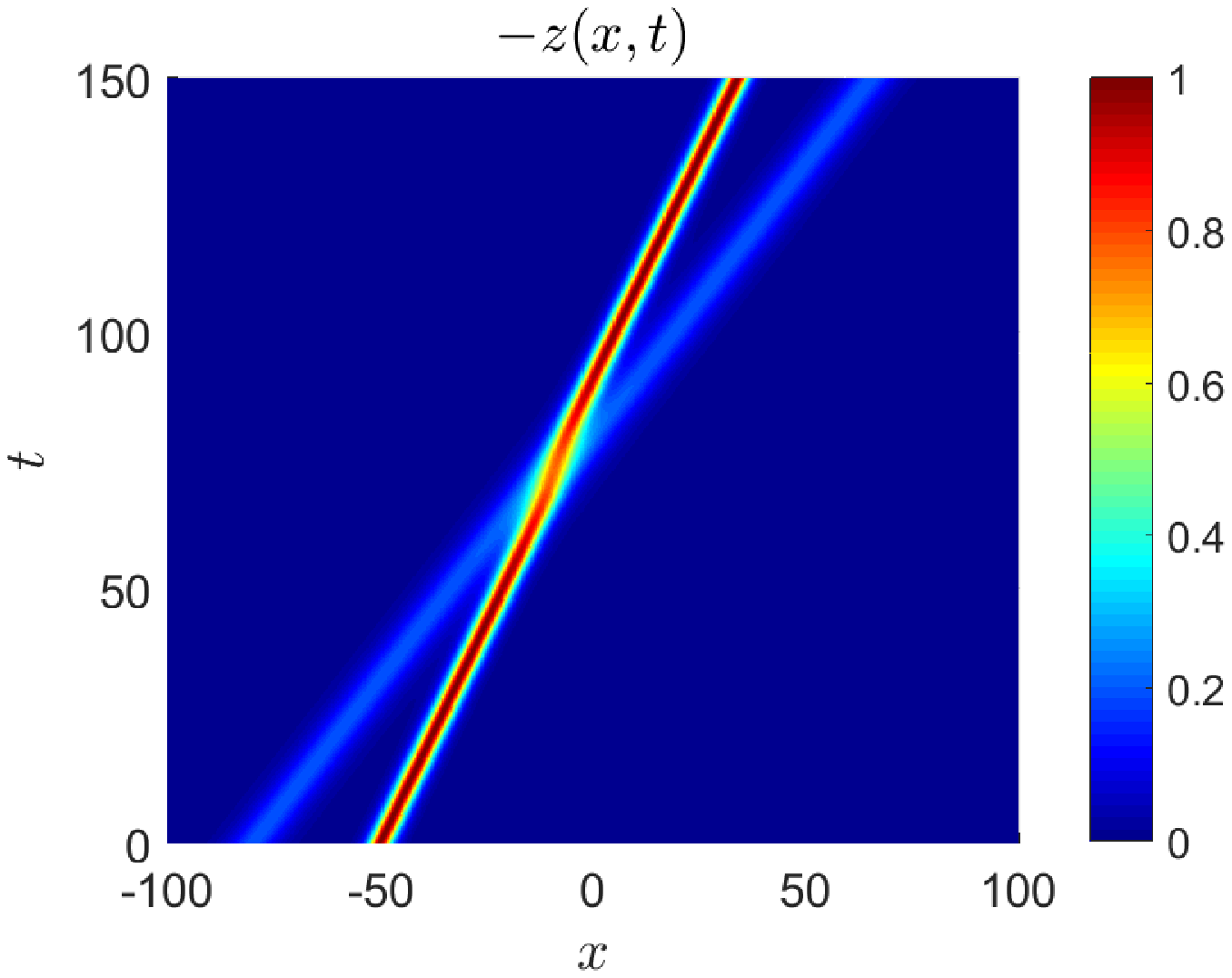}
\end{minipage}%
\hspace{3mm}
\begin{minipage}[t]{0.5\linewidth}
\centering
\includegraphics[height=4cm,width=7.1cm]{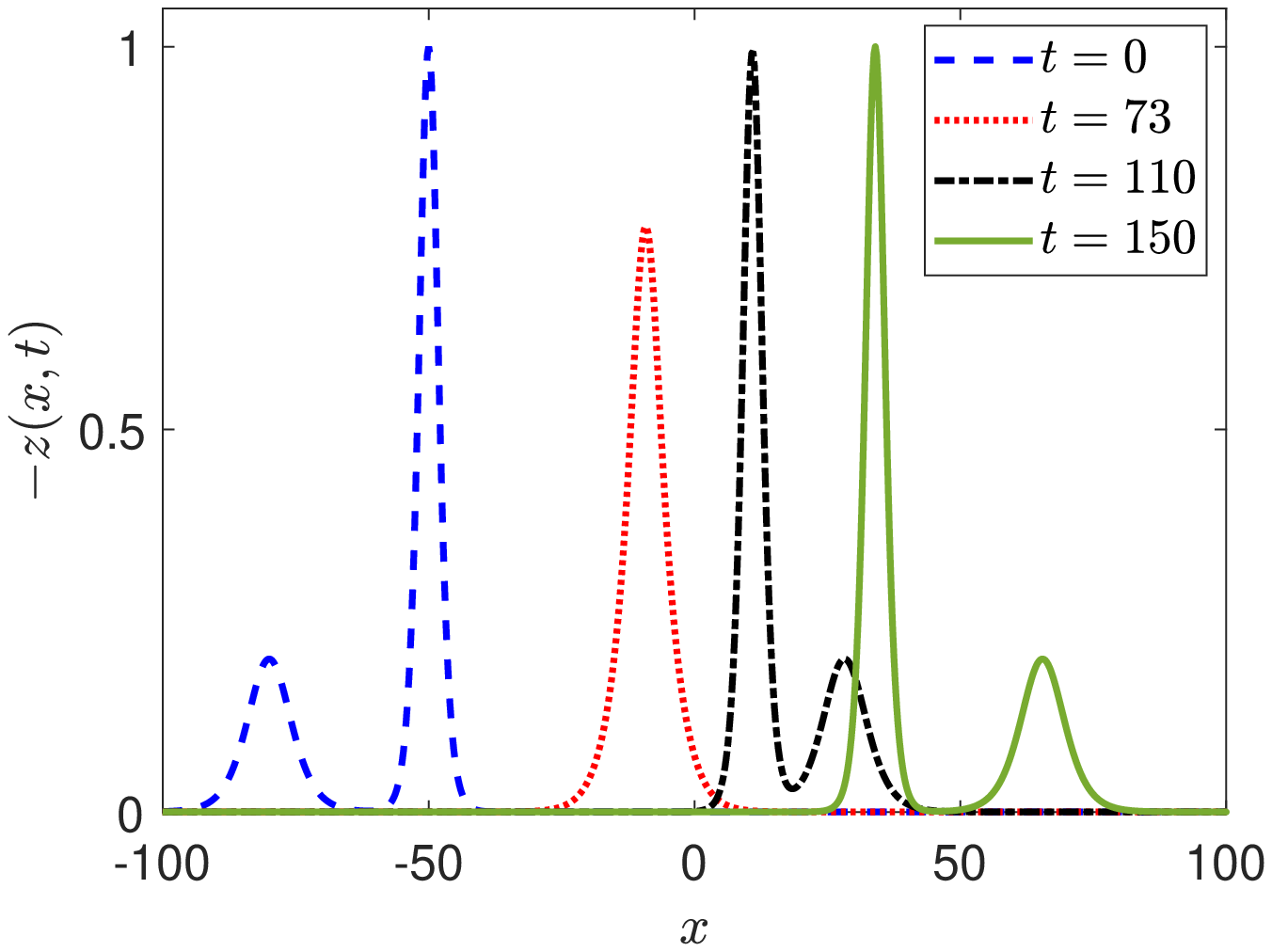}
\end{minipage}
\vspace{-3mm}
\caption{Overtaking interaction of two solitons traveling in the same direction.}\label{chase}
\end{figure}

Fig. \ref{chase} displays the \emph{overtaking interaction} of two solitons moving in the same direction with different velocities. The faster wave catches up with the slower one and leaves it behind as time evolves. Both solitons recover their shape and velocities after interaction. Noticing that the amplitude decreases during the interaction, which is completely different from the \emph{head-on interaction} case, where the amplitude is strengthened during the collision (cf. Fig. \ref{interac}).

\section{Conclusions}
A Deuflhard-type exponential integrator was proposed and analyzed for the ``Good" Boussinesq (GB) equation with a general nonlinearity. The method was proved to unconditionally converge at the second order in time and spectrally in space, respectively, in a generic energy norm. Specifically, it requires four additional orders of regularities on the solution to attain the quadratic convergence rate in time. Numerical results confirm our analytical results. Extensive numerical experiments reveal rich and complicated dynamical phenomena for the GB equation, such as \emph{elastic interaction}, \emph{blow-up phenomena}, \emph{recovering as a static soliton} and so forth.


\begin{thebibliography}{100}
\bibitem{Adams}
Adams R. A., Fournier J. J.: Sobolev spaces. Elsevier (2003)


\bibitem{Bona1988}
Bona J.  L., Sachs R. L.: Global existence of smooth solutions and stability of solitary
waves for a generalized Boussinesq equation. Comm. Math. Phys. \textbf{118}, 15-29 (1988)

\bibitem{Bona}
Bona J. L., Smith R. A.: A model for the two-way propagation of water waves in a channel. Math. Proc. Cambridge Philos. Soc. \textbf{79}, 167-182 (1976)

\bibitem{boussinesq1872}
Boussinesq J.: Th{\'e}orie des ondes et des remous qui se propagent le long d'un canal rectangulaire horizontal, en communiquant au liquide contenu dans ce canal des vitesses sensiblement pareilles de la surface au fond. J. Math. Pures Appl. \textbf{17}, 55-108 (1872)

\bibitem{bratsos2007}
 Bratsos A. G.: A second order numerical scheme for the solution of the one-dimensional Boussinesq equation. Numer. Algor. \textbf{46}, 45-58 (2007)

 \bibitem{cai2013}
 Cai J., Wang Y.: Local structure-preserving algorithms for the ``good" Boussinesq equation. J. Comp. Phys. \textbf{239}, 72-89 (2013)

\bibitem{Super}
Chartier Ph., M\'{e}hats F., Thalhammer M., Zhang Y.: Improved error estimates for splitting methods applied to highly-oscillatory nonlinear Schr\"{o}dinger equations. Math. Comp. \textbf{85}, 2863-2885 (2016)

\bibitem{Chen2017}
Chen M, Kong L., Hong Y.: Efficient structure-preserving schemes for good Boussinesq equation. Math. Meth. Appl. Sci. \textbf{41}, 1743-1752 (2018)

\bibitem{cheng2015}
Cheng K., Feng W., Gottlieb S., Wang C.: A Fourier pseudospectral method for the ``good" Boussinesq equation with second-order temporal accuracy. Numer. Methods Partial Differential Equations \textbf{31}, 202-224 (2015)

\bibitem{dehghan2012}
Dehghan M., Salehi R.: A meshless based numerical technique for traveling solitary wave solution of Boussinesq equation. Appl. Math. Model. \textbf{36}, 1939-1956 (2012)

\bibitem{Deuf}
Deuflhard P.: A study of extrapolation methods based on multistep schemes without parasitic solutions. ZAMP \textbf{30}, 177-189 (1979)

\bibitem{Zoh}
El-Zoheiry H.: Numerical investigation for the solitary waves interaction of the
``good" Boussinesq equation. Appl. Numer. Math. \textbf{45}, 161-173 (2003)

\bibitem{fang1996}
Fang Y., Grillakis M.: Existence and uniqueness for Boussinesq type equations on a circle. Comm. Partial Differential Equations \textbf{21}, 1253-1277 (1996)

\bibitem{farah2009}
Farah L.:  Local solutions in Sobolev spaces with negative indices for the ``good" Boussinesq equation, Comm. Partial Differential Equations \textbf{34}, 52-73 (2009)

\bibitem{farah2010}
Farah L., Scialom M.: On the periodic ``good" Boussinesq equation. Proc. Amer. Math. Soc. \textbf{138}, 953-964 (2010)

\bibitem{de1990}
Frutos J.~De, Ortega T., Sanz-Serna J.~M.: A Hamiltonian explicit algorithm with
spectral accuracy for the ``good" Boussinesq equation. Comput. Methods Appl.
Mech. Engrg. \textbf{80}, 417-423 (1990)

\bibitem{de1991}
Frutos J.~De, Ortega T., Sanz-Serna J.~M.: Pseudospectral method
  for the ``good" Boussinesq equation. Math. Comput. \textbf{57}, 109-122 (1991)

\bibitem{exp10}
Hochbruck M., Ostermann A.: Exponential integrators, Acta Numer. 209-286 (2010)

\bibitem{Ismail}
Ismail M. S., Mosally F.: A fourth order finite difference method for the good Boussinesq
equation. Abstr. Appl. Anal. (2014)

\bibitem{jiang2016}
Jiang C., Sun J., He X., Zhou L.: High order energy-preserving method of the ``good" Boussinesq equation. Numer. Math. Theor. Meth. Appl. \textbf{9}, 111-122 (2016)

\bibitem{kishimoto2013sharp}
Kishimoto N.: Sharp local well-posedness for the ``good"
Boussinesq equation. J. Differential Equations \textbf{254}, 2393-2433 (2013)

\bibitem{Kish2010}
Kishimoto N., Tsugawa K.: Local well-posedness for quadratic nonlinear Schr\"odinger equations and the ``good" Boussinesq equation. Differential Integral Equations \textbf{23}, 463-493 (2010)

\bibitem{manoranjan1984}
Manoranjan V.~S., Mitchell A., Morris J.~L.: Numerical solutions of the good Boussinesq equation. SIAM J. Sci. Comput. \textbf{5}, 946-957 (1984)

\bibitem{manoranjan1988}
Manoranjan V.~S., Ortega T., Sanz-Serna J.~M.: Soliton and antisoliton interactions in the ``good" Boussinesq equation. J. Math. Phys. \textbf{29}, 964-1968 (1988)

\bibitem{Akbar}
Mohebbi A.,  Asgari Z.: Efficient numerical algorithms for the solution of ``good" Boussinesq equation in water wave propagation. Comp. Phys. Comm. \textbf{182}, 2464-2470 (2011)

\bibitem{oh2013improved}
Oh S., Stefanov A.: Improved local well-posedness for the
  periodic ``good" Boussinesq equation. J. Differential Equations \textbf{254}, 4047-4065  (2013)

\bibitem{ortega1990}
Ortega T., Sanz-Serna J.~M.: Nonlinear stability and convergence   of finite-difference methods for the ``good" Boussinesq equation. Numer. Math. \textbf{58}, 215-229 (1990)

\bibitem{Alex19}
Ostermann A., Su C.: Two exponential-type integrators for the ``good" Boussinesq equation. Numer. Math. doi.org/10.1007/s00211-019-01064-4 (2019)

\bibitem{shen}
Shen J., Tang T.: Spectral and high-order methods with applications. Science Press, Beijing (2006)

\bibitem{Varl}
Varlamov V.: Eigenfunction expansion method and the long-time asymptotics for the damped Boussinesq equation. Discrete Contin. Dyn. Syst. \textbf{7}, 675-702 (2001)

\bibitem{yan2016}
Yan J., Zhang Z.: New energy-preserving schemes using Hamiltonian boundary value and Fourier pseudospectral methods for the   numerical solution of the ``good" Boussinesq equation. Comput. Phys. Commun. \textbf{201}, 33-42 (2016)

\bibitem{zhang2018}
Zhang C., Huang J., Wang C., Yue X.: On the operator splitting
  and integral equation preconditioned deferred correction methods for the
  ``Good" Boussinesq equation. J. Sci. Comput. \textbf{75}, 687-712 (2018)

\bibitem{zhang2017}
Zhang C., Wang H., Huang J., Wang C., Yue X.: A second order
  operator splitting numerical scheme for the ``good" Boussinesq equation. Appl. Numer. Math. \textbf{119}, 179-193 (2017)

\bibitem{zhao2016error}
Zhao X.: On error estimates of an exponential wave integrator sine pseudospectral method for the Klein-Gordon-Zakharov system. Numer. Methods Partial Differential Equations \textbf{32}, 266-291 (2016)

\end{thebibliography}
\end{document}